\documentclass[leqno,12pt]{article}

\usepackage{pstricks}
\usepackage{amssymb}
\usepackage{mathrsfs} 
\usepackage{amsmath}
\usepackage{amsthm}
\usepackage{amsfonts}
\usepackage{hyperref} 
\usepackage{verbatim} 
\usepackage{graphicx}

\usepackage{pgf,tikz}
\usetikzlibrary{arrows}
\usepackage{caption}
\usetikzlibrary{decorations}
\usetikzlibrary{decorations.pathreplacing}
\tikzset{individu/.style={fill,thick,circle}}
\usepackage{pgfplots}

\def\thesection{\arabic{section}}
\renewcommand{\theequation}{\thesection.\arabic{equation}}

\newtheorem{theorem}{Theorem}[section]
\newtheorem{lemma}[theorem]{Lemma}
\newtheorem{proposition}[theorem]{Proposition}

\newtheorem{fact}[theorem]{Fact}

\newtheorem{definition}[theorem]{Definition}

\theoremstyle{definition}   
\newtheorem{remark}[theorem]{Remark}
\newtheorem{example}[theorem]{Example}

\evensidemargin\oddsidemargin

\newcommand{\eqnsection}{
\renewcommand{\theequation}{\thesection.\arabic{equation}}
    \makeatletter
    \csname  @addtoreset\endcsname{equation}{section}
    \makeatother}
\eqnsection

\def\r{{\mathbb R}}

\def\p{{\mathbb P}}

\def\P{{\bf P}}
\def\E{{\bf E}}

\def\z{{\mathbb Z}}

\def\ee{\mathrm{e}}

\def\T{{\mathbb T}}

\begin{document}

\baselineskip=18pt
\setcounter{page}{1}


\vglue50pt

\centerline{\large\bf The sustainability probability for the critical Derrida--Retaux model}



\bigskip
\bigskip

\centerline{Xinxing Chen,\footnote{\scriptsize School of Mathematical Sciences, Shanghai Jiaotong University, 200240 Shanghai, China, \texttt{chenxinx@sjtu.edu.cn} $\;$Partially supported by NSFC grants 11771286 and 11531001.}$\;$
Yueyun Hu,\footnote{\scriptsize LAGA, Universit\'e Sorbonne Paris Nord, 99 av.~J-B Cl\'ement, F-93430 Villetaneuse, France, \texttt{yueyun@math.univ-paris13.fr} $\;$Partially supported by ANR MALIN and ANR SWIWS.}$\;$
and Zhan Shi\footnote{\scriptsize LPSM, Sorbonne Universit\'e Paris VI, 4 place Jussieu, F-75252 Paris Cedex 05, France, \texttt{zhan.shi@sorbonne-universite.fr} $\;$Partially supported by ANR MALIN.}
} 



\bigskip
\bigskip

{\leftskip=1.5truecm \rightskip=1.5truecm \baselineskip=15pt \small

\noindent{\slshape\bfseries Summary.} We are interested in the recursive model $(Y_n, \, n\ge 0)$ studied by Collet, Eckmann, Glaser and Martin~\cite{collet-eckmann-glaser-martin} and by Derrida and Retaux~\cite{derrida-retaux}. We prove that at criticality, the probability $\P(Y_n>0)$ behaves like $n^{-2 + o(1)}$ as $n$ goes to infinity; this gives a weaker confirmation of predictions made in \cite{collet-eckmann-glaser-martin}, \cite{derrida-retaux} and \cite{bmxyz_questions}. Our method relies on studying the number of pivotal vertices and open paths, combined with a delicate coupling argument.

\bigskip


\bigskip


} 

\bigskip
\bigskip

\section{Introduction}
\label{s:intro}

Let us start with a simple example.

\medskip

\begin{example}
 \label{ex:DR}

 Let $Y_0$ be a random variable with $\P (Y_0=2) = \frac15 = 1- \P (Y_0=0)$. Consider the recurrence relation as follows: for all $n\ge 0$,
 $$
     Y_{n+1}
     =
     (Y_n + \widehat{Y}_n -1)^+ ,
 $$

 \noindent where $\widehat{Y}_n$ is an independent copy of $Y_n$, and $x^+ := \max\{ x, \, 0\}$ for all $x\in \r$. Then, from the physics literature, one expects
 \begin{eqnarray}
     \P(Y_n >0) 
  &\sim& \frac{4}{n^2},
     \label{4/n2}
     \\
     \E(Y_n) 
  &\sim& \frac{8}{n^2}, \qquad n\to \infty\, .
     \label{8/n2}
 \end{eqnarray}
 
 \noindent [Notation: $a_n \sim b_n$ means $\lim_{n\to \infty} \frac{a_n}{b_n} =1$.] See \cite{collet-eckmann-glaser-martin}, \cite{derrida-retaux} and \cite{bmxyz_questions} for explanations of \eqref{4/n2}, and \cite{bmxyz_questions} for explanations of \eqref{8/n2}.
 
 However, from a mathematical point of view, \eqref{4/n2} and \eqref{8/n2} are mere conjectures. The only relevant rigorous result (see \cite{bmxyz_questions}) available is:
 $$
 \limsup_{n\to \infty} \, n \, \P(Y_n >0) < \infty \, .
 $$
 
 \noindent The main result of our paper, stated in Theorem \ref{t:main} below, will imply that
 $$
 \P(Y_n >0) = \frac{1}{n^{2+o(1)}},
 \qquad
 \E(Y_n) = \frac{1}{n^{2+o(1)}},
 \qquad
 n\to \infty \, .
 $$
 
 \noindent We are going to see that this is valid for a large family of initial distributions of $Y_0$.\qed

\end{example}

\medskip

Throughout the paper, we fix an integer $m\ge 2$, and let $X_0 \ge 0$ be a random variable taking values in $\z_+ := \{ 0, \, 1, \, 2, \ldots\}$ and satisfying 
$\E(X_0 \, m^{X_0})<\infty$. We consider the following recurrence relation: for all $n\ge 0$,
\begin{equation}
    X_{n+1}
    =
    (X_n^{(1)} + \cdots + X_n^{(m)} -1)^+ ,
    \label{iteration}
\end{equation}

\noindent where $X_n^{(i)}$, $i\ge 1$, are independent copies of $X_n$. The assumption $\E(X_0 \, m^{X_0})<\infty$ ensures that $\E(X_n \, m^{X_n})<\infty$ for all $n\ge 0$. The basic question is how the distribution of $X_n$ behaves when $n$ is sufficiently large. 

The random recursive system $(X_n, \, n\ge 0)$, referred to as the Derrida--Retaux system or the Derrida--Retaux model, has been investigated by Derrida and Retaux~\cite{derrida-retaux} as a simplified hierarchical renormalization model in order to understand the depinning transition of a line in presence of strong disorder. Pinning models have an intensive literature both in mathematics and in theoretical physics; see for example \cite{alexander}, \cite{berger-giacomin-lacoin}, \cite{giacomin}, \cite{giacomin_stf}, \cite{giacomin-toninelli}, \cite{giacomin-toninelli-lacoin}, and \cite{derrida-hakim-vannimenus}, \cite{derrida-retaux}, \cite{monthus}, \cite{tang-charte}. We refer to \cite{derrida-retaux} and \cite{bz_DRsurvey} for more references and explanations of the recursion \eqref{iteration}. The Derrida--Retaux system has also appeared in other contexts: in Collet et al.~\cite{collet-eckmann-glaser-martin} as a spin glass model, in Li and Rogers~\cite{li-rogers} as an iteration function of random variables, in Goldschmidt and Przykucki~\cite{goldschmidt-przykucki} and Curien and H\'enard~\cite{curien-henard} as a parking scheme; it also belongs to one of the max-type recursion families in the seminal paper by Aldous and Bandyopadhyay~\cite{aldous-bandyopadhyay}. 

The system is known for exhibiting a phase transition:\footnote{It is crucial to assume that in the Derrida--Retaux model, $X_0$ is integer-valued; it is an open problem to extend Theorem A without this assumption. See Derrida and Retaux~\cite{derrida-retaux} for more discussions.}

\bigskip

\noindent {\bf Theorem A. (Collet et al.~\cite{collet-eckmann-glaser-martin})} {\it Assume $\E(X_0 \, m^{X_0})<\infty$. 

{\rm (i)} If $(m-1) \E(X_0 \, m^{X_0}) >\E(m^{X_0})$, then $\lim_{n\to \infty} \, \frac{\E(X_n)}{m^n}$ exists and is positive.

{\rm (ii)} If $(m-1) \E(X_0 \, m^{X_0}) \le \E(m^{X_0})$, then $\sup_{n\ge 0} \E(X_n) \le \frac{1}{m-1}$.}

\bigskip

According to Theorem A, there is a dichotomy for the asymptotic behaviour of $\E(X_n)$: it either grows exponentially fast as $n\to \infty$, or is always bounded by the constant $\frac{1}{m-1}$. The system $(X_n, \, n\ge 0)$ is said to be supercritical if $(m-1) \E(X_0 \, m^{X_0}) >\E(m^{X_0})$, critical if $(m-1) \E(X_0 \, m^{X_0}) = \E(m^{X_0})<\infty$, and subcritical if $(m-1) \E(X_0 \, m^{X_0}) <\E(m^{X_0})$.

The present paper focuses on the critical case: $(m-1) \E(X_0 \, m^{X_0}) = \E(m^{X_0})<\infty$. In order to insist on the criticality, we write $(Y_n, \, n\ge 0)$ instead of $(X_n, \, n\ge 0)$ for the critical Derrida--Retaux system. It is known (see \cite{collet-eckmann-glaser-martin}, \cite{bmxyz_questions}) that
$$
\P(Y_n \ge 1) \to 0, \qquad n\to \infty \, .
$$

\noindent We study the rate of decay of $\P(Y_n \ge 1)$.

\medskip

\begin{theorem}
 \label{t:main} 
 
 If $\E(t^{Y_0})<\infty$ for some $t>m$, then for $n\to \infty$,
 $$
 \P(Y_n\ge 1)
 = 
 \frac{1}{n^{2+o(1)}},
 \quad 
 \E(Y_n)
 = 
 \frac{1}{n^{2+o(1)}} \, .
 $$

\end{theorem}

\medskip

More precisely, our Proposition \ref{p:ub} implies that if $\E(Y_0^5 m^{Y_0})<\infty$, then $\E(Y_n)$ is bounded from above by a constant multiple of $\frac{(\log n)^2}{n^2}$ for all $n\ge 2$, whereas Proposition \ref{p:lb} yields that if $\E(t^{Y_0})<\infty$ for some $t>m$, then $\P(Y_n\ge 1)$ is bounded from below by a constant multiple of $\frac{1}{n^2(\log n)^\nu}$ for some $\nu>0$ and all $n\ge 2$.

The problem of evaluating $\P(Y_n\ge 1)$ was discussed in \cite{collet-eckmann-glaser-martin}, \cite{derrida-retaux} and \cite{bmxyz_questions}, and $\E(Y_n)$ in \cite{bmxyz_questions}; it was expected, without precision on integrability condition on $Y_0$, that
\begin{equation}
    \P(Y_n\ge 1) \sim \frac{c_*}{n^2}, 
    \qquad
    \E(Y_n) \sim \frac{c_{**}}{n^2}, 
    \qquad
    n\to \infty \, ,
    \label{prediction_P(Y>0)}
\end{equation}

\noindent with $c_* := \frac{4}{(m-1)^2}$ and $c_{**} := \frac{m}{m-1} \, c_*$. In \cite{bmvxyz_conjecture_DR}, it was realized that \eqref{prediction_P(Y>0)} should be valid only under the assumption $\E(Y_0^3 m^{Y_0})<\infty$. When the latter condition fails, $\P(Y_n\ge 1)$ and $\E(Y_n)$ should behave differently. For example, for $m=2$, if $\P(Y_0=k) \sim c_0 \, k^{-\alpha} 2^{-k}$, $k\to \infty$, for some $\alpha \in (2, \, 4]$ (so $\E(Y_0^3 m^{Y_0}) =\infty$) and $c_0 \in (0, \, \infty)$, then it is expected (see \cite{bz_DRsurvey}) that
\begin{equation}
    \P(Y_n\ge 1) \sim \frac{c(\alpha)}{n^2}, \qquad n\to \infty \,
    \label{alpha}
\end{equation}

\noindent with $c(\alpha) := \frac{\alpha(\alpha-2)}{2}$. For discussions on other related questions, see \cite{bz_DRsurvey}. 

Let us also mention that some continuous-time systems were studied in \cite{derrida-retaux}, \cite{HMP} and \cite{4authors}, for which the analogue of the asymptotic equivalences of $\P(Y_n\ge 1)$ was proved; see \cite{derrida-retaux} and \cite{HMP} for the analogue of \eqref{prediction_P(Y>0)}, and \cite{4authors} for the analogue of \eqref{alpha}.

When a Derrida--Retaux system $(X_n, \, n\ge 0)$ is subcritical, it is not hard to prove that $\P(X_n \ge 1)$ and $\E(X_n)$ decay exponentially (as $n\to \infty$) if $X_0$ satisfies an appropriate integrability condition.


The proof of Theorem \ref{t:main} relies on the study of the important notions of {\it pivotal vertices} and {\it open paths} which we define in Section \ref{s:pivot}; these notions are closely related. The main advantage of studying the expected number of pivotal vertices is that it satisfies a nice identity. This is Theorem \ref{t:new_iteration}, valid for any Derrida--Retaux system (regardless of whether it is supercritical, critical or subcritical) and useful in the proof of both upper and lower bounds in Theorem \ref{t:main}. In the proof of the upper bound in Theorem \ref{t:main}, only a simple consequence (inequality \eqref{E(N_Y>k)}) of Theorem \ref{t:new_iteration} is needed a couple of times. The proof of the lower bound relies on the full strength of Theorem \ref{t:new_iteration} and goes actually beyond, in the sense that pivotal vertices are studied in a more complicated setting, and that the key step in the proof consists in establishing a weaker form of inequality \eqref{E(N_Y>k)} in the opposite direction.

The notion of open paths was already introduced in \cite{bmvxyz_conjecture_DR}, and was known to be an important tool in the study of the Derrida--Retaux model. {F}rom methodological point of view, an interesting contribution of the paper is to introduce the notion of pivotal vertices, and to prove Theorem \ref{t:new_iteration} on the expected number of pivotal vertices. Another notable contribution of the paper consists in studying the Derrida--Retaux system at intermediate generations. The idea of studying intermediate generations is quite natural in statistical physics; our in-depth descriptions (in Sections \ref{s:ub} and \ref{s:lb}) of Derrida--Retaux systems at appropriate intermediate generations are among the key ingredients in the proof of Theorem \ref{t:main}. These descriptions require a careful analysis of open paths and pivotal vertices and of the moment generating function.

The rest of the paper is as follows. In Section \ref{s:pivot}, we introduce pivotal vertices and open paths, and prove the aforementioned Theorem \ref{t:new_iteration}. Section \ref{s:preparation_ub} collects some preliminary estimates for the moment generating function of critical Derrida--Retaux systems; these estimates are either known, or are routinely proved by means of existing arguments. The upper and lower bounds in Theorem \ref{t:main} are proved in Sections \ref{s:ub} and \ref{s:lb}, respectively.

\section{Pivotal vertices and open paths}
\label{s:pivot}

Let $(X_n, \, n\ge 0)$ be a Derrida--Retaux system defined via \eqref{iteration}. 
There is a natural hierarchical representation of the system (see \cite{collet-eckmann-glaser-martin2}, \cite{derrida-retaux}, \cite{yz_bnyz}).


Let us consider a family of random variables $(X(v), \, v\in \T)$, indexed by a (reversed) $m$-ary tree $\T$, defined as follows. For any vertex $v$ in the genealogical tree $\T$, let $|v|$ denote the generation of $v$ (so $|v|=0$ if the vertex $v$ is in the initial generation). Let $X(v)$, for $v\in \T$ with $|v|=0$, be independent and identically distributed (i.i.d.) having the law of $X_0$. For any $v\in \T$ with $|v|\ge 1$, let $v^{(1)}$, $\ldots$, $v^{(m)}$ denote the $m$ parents of $v$ in generation $|v|-1$, and define
$$
X(v)
:=
(X(v^{(1)}) + \cdots + X(v^{(m)})-1)^+ \, .
$$

\noindent By definition, for any $n\ge 0$, the random variables $X(v)$, for $v\in \T$ with $|v|=n$, are i.i.d.\ having the law of $X_n$.

For $n\ge 0$, let $\mathfrak{e}_n$ denote the first lexicographic vertex in the $n$-th generation of $\T$. Let $\T_n$ denote the (reversed) subtree formed by all the ancestors (including $\mathfrak{e}_n$ itself) of $\mathfrak{e}_n$ in the first $n$ generations. See Figure \ref{f:fig1} for an example.

\bigskip

\begin{figure}[ht]
\centering
\begin{tikzpicture}[scale=1]

    \draw (0,8) -- (0,7.8) ; 
    \draw (-0.2,8) -- (0.2,8) -- (0.2,8.4) -- (-0.2,8.4) -- cycle;
    \draw (0,7.94) node[above] {\footnotesize 1} ;
    \draw (0.8,8) -- (0.8,7.8) ; 
    \draw (0.8-0.2,8) -- (0.8+0.2,8) -- (0.8+0.2,8.4) -- (0.8-0.2,8.4) -- cycle;
    \draw (0.8,7.94) node[above] {\footnotesize 0} ;
    \draw (0,7.8) -- (0.8,7.8) ; 
    \draw (0.4,7.8) -- (0.4,7.6) ; 
    \draw (0.4-0.2,7.6) -- (0.4+0.2,7.6) -- (0.4+0.2,7.2) -- (0.4-0.2,7.2) -- cycle;
    \draw (0.4,7.14) node[above] {\footnotesize 0} ;

    \draw (1.6,8) -- (1.6,7.8) ; 
    \draw (1.6-0.2,8) -- (1.6+0.2,8) -- (1.6+0.2,8.4) -- (1.6-0.2,8.4) -- cycle;
    \draw (1.6,7.94) node[above] {\footnotesize 0} ;
    \draw (1.6+0.8,8) -- (1.6+0.8,7.8) ; 
    \draw (1.6+0.8-0.2,8) -- (1.6+0.8+0.2,8) -- (1.6+0.8+0.2,8.4) -- (1.6+0.8-0.2,8.4) -- cycle;
    \draw (1.6+0.8,7.94) node[above] {\footnotesize 0} ;
    \draw (1.6,7.8) -- (1.6+0.8,7.8) ; 
    \draw (1.6+0.4,7.8) -- (1.6+0.4,7.6) ; 
    \draw (1.6+0.4-0.2,7.6) -- (1.6+0.4+0.2,7.6) -- (1.6+0.4+0.2,7.2) -- (1.6+0.4-0.2,7.2) -- cycle;
    \draw (1.6+0.4,7.14) node[above] {\footnotesize 0} ;

    \draw (3.2,8) -- (3.2,7.8) ; 
    \draw (3.2-0.2,8) -- (3.2+0.2,8) -- (3.2+0.2,8.4) -- (3.2-0.2,8.4) -- cycle;
    \draw (3.2,7.94) node[above] {\footnotesize 3} ;
    \draw (3.2+0.8,8) -- (3.2+0.8,7.8) ; 
    \draw (3.2+0.8-0.2,8) -- (3.2+0.8+0.2,8) -- (3.2+0.8+0.2,8.4) -- (3.2+0.8-0.2,8.4) -- cycle;
    \draw (3.2+0.8,7.94) node[above] {\footnotesize 3} ;
    \draw (3.2,7.8) -- (3.2+0.8,7.8) ; 
    \draw (3.2+0.4,7.8) -- (3.2+0.4,7.6) ; 
    \draw (3.2+0.4-0.2,7.6) -- (3.2+0.4+0.2,7.6) -- (3.2+0.4+0.2,7.2) -- (3.2+0.4-0.2,7.2) -- cycle;
    \draw (3.2+0.4,7.14) node[above] {\footnotesize 5} ;

    \draw (4.8,8) -- (4.8,7.8) ; 
    \draw (4.8-0.2,8) -- (4.8+0.2,8) -- (4.8+0.2,8.4) -- (4.8-0.2,8.4) -- cycle;
    \draw (4.8,7.94) node[above] {\footnotesize 0} ;
    \draw (4.8+0.8,8) -- (4.8+0.8,7.8) ; 
    \draw (4.8+0.8-0.2,8) -- (4.8+0.8+0.2,8) -- (4.8+0.8+0.2,8.4) -- (4.8+0.8-0.2,8.4) -- cycle;
    \draw (4.8+0.8,7.94) node[above] {\footnotesize 1} ;
    \draw (4.8,7.8) -- (4.8+0.8,7.8) ; 
    \draw (4.8+0.4,7.8) -- (4.8+0.4,7.6) ; 
    \draw (4.8+0.4-0.2,7.6) -- (4.8+0.4+0.2,7.6) -- (4.8+0.4+0.2,7.2) -- (4.8+0.4-0.2,7.2) -- cycle;
    \draw (4.8+0.4,7.14) node[above] {\footnotesize 0} ;

    \draw (6.4,8) -- (6.4,7.8) ; 
    \draw (6.4-0.2,8) -- (6.4+0.2,8) -- (6.4+0.2,8.4) -- (6.4-0.2,8.4) -- cycle;
    \draw (6.4,7.94) node[above] {\footnotesize 0} ;
    \draw (6.4+0.8,8) -- (6.4+0.8,7.8) ; 
    \draw (6.4+0.8-0.2,8) -- (6.4+0.8+0.2,8) -- (6.4+0.8+0.2,8.4) -- (6.4+0.8-0.2,8.4) -- cycle;
    \draw (6.4+0.8,7.94) node[above] {\footnotesize 0} ;
    \draw (6.4,7.8) -- (6.4+0.8,7.8) ; 
    \draw (6.4+0.4,7.8) -- (6.4+0.4,7.6) ; 
    \draw (6.4+0.4-0.2,7.6) -- (6.4+0.4+0.2,7.6) -- (6.4+0.4+0.2,7.2) -- (6.4+0.4-0.2,7.2) -- cycle;
    \draw (6.4+0.4,7.14) node[above] {\footnotesize 0} ;

    \draw (8.0,8) -- (8.0,7.8) ; 
    \draw (8.0-0.2,8) -- (8.0+0.2,8) -- (8.0+0.2,8.4) -- (8.0-0.2,8.4) -- cycle;
    \draw (8.0,7.94) node[above] {\footnotesize 0} ;
    \draw (8.0+0.8,8) -- (8.0+0.8,7.8) ; 
    \draw (8.0+0.8-0.2,8) -- (8.0+0.8+0.2,8) -- (8.0+0.8+0.2,8.4) -- (8.0+0.8-0.2,8.4) -- cycle;
    \draw (8.0+0.8,7.94) node[above] {\footnotesize 3} ;
    \draw (8.0,7.8) -- (8.0+0.8,7.8) ; 
    \draw (8.0+0.4,7.8) -- (8.0+0.4,7.6) ; 
    \draw (8.0+0.4-0.2,7.6) -- (8.0+0.4+0.2,7.6) -- (8.0+0.4+0.2,7.2) -- (8.0+0.4-0.2,7.2) -- cycle;
    \draw (8.0+0.4,7.14) node[above] {\footnotesize 2} ;

    \draw (9.6,8) -- (9.6,7.8) ; 
    \draw (9.6-0.2,8) -- (9.6+0.2,8) -- (9.6+0.2,8.4) -- (9.6-0.2,8.4) -- cycle;
    \draw (9.6,7.94) node[above] {\footnotesize 0} ;
    \draw (9.6+0.8,8) -- (9.6+0.8,7.8) ; 
    \draw (9.6+0.8-0.2,8) -- (9.6+0.8+0.2,8) -- (9.6+0.8+0.2,8.4) -- (9.6+0.8-0.2,8.4) -- cycle;
    \draw (9.6+0.8,7.94) node[above] {\footnotesize 0} ;
    \draw (9.6,7.8) -- (9.6+0.8,7.8) ; 
    \draw (9.6+0.4,7.8) -- (9.6+0.4,7.6) ; 
    \draw (9.6+0.4-0.2,7.6) -- (9.6+0.4+0.2,7.6) -- (9.6+0.4+0.2,7.2) -- (9.6+0.4-0.2,7.2) -- cycle;
    \draw (9.6+0.4,7.14) node[above] {\footnotesize 0} ;

    \draw (11.2,8) -- (11.2,7.8) ; 
    \draw (11.2-0.2,8) -- (11.2+0.2,8) -- (11.2+0.2,8.4) -- (11.2-0.2,8.4) -- cycle;
    \draw (11.2,7.94) node[above] {\footnotesize 1} ;
    \draw (11.2+0.8,8) -- (11.2+0.8,7.8) ; 
    \draw (11.2+0.8-0.2,8) -- (11.2+0.8+0.2,8) -- (11.2+0.8+0.2,8.4) -- (11.2+0.8-0.2,8.4) -- cycle;
    \draw (11.2+0.8,7.94) node[above] {\footnotesize 0} ;
    \draw (11.2,7.8) -- (11.2+0.8,7.8) ; 
    \draw (11.2+0.4,7.8) -- (11.2+0.4,7.6) ; 
    \draw (11.2+0.4-0.2,7.6) -- (11.2+0.4+0.2,7.6) -- (11.2+0.4+0.2,7.2) -- (11.2+0.4-0.2,7.2) -- cycle;
    \draw (11.2+0.4,7.14) node[above] {\footnotesize 0} ;

    \draw (0.4,7.2) -- (0.4,7.0) ; 
    \draw (2,7.2) -- (2,7.0) ; 
    \draw (0.4,7.0) -- (2,7.0) ; 
    \draw (1.2,7.0) -- (1.2,6.8) ; 
    \draw (1.2-0.2,6.8) -- (1.2+0.2,6.8) -- (1.2+0.2,6.4) -- (1.2-0.2,6.4) -- cycle;
    \draw (1.2,6.35) node[above] {\footnotesize 0} ;

    \draw (3.2+0.4,7.2) -- (3.2+0.4,7.0) ; 
    \draw (3.2+2,7.2) -- (3.2+2,7.0) ; 
    \draw (3.2+0.4,7.0) -- (3.2+2,7.0) ; 
    \draw (3.2+1.2,7.0) -- (3.2+1.2,6.8) ; 
    \draw (3.2+1.2-0.2,6.8) -- (3.2+1.2+0.2,6.8) -- (3.2+1.2+0.2,6.4) -- (3.2+1.2-0.2,6.4) -- cycle;
    \draw (3.2+1.2,6.35) node[above] {\footnotesize 4} ;

    \draw (6.4+0.4,7.2) -- (6.4+0.4,7.0) ; 
    \draw (6.4+2,7.2) -- (6.4+2,7.0) ; 
    \draw (6.4+0.4,7.0) -- (6.4+2,7.0) ; 
    \draw (6.4+1.2,7.0) -- (6.4+1.2,6.8) ; 
    \draw (6.4+1.2-0.2,6.8) -- (6.4+1.2+0.2,6.8) -- (6.4+1.2+0.2,6.4) -- (6.4+1.2-0.2,6.4) -- cycle;
    \draw (6.4+1.2,6.35) node[above] {\footnotesize 1} ;

    \draw (9.6+0.4,7.2) -- (9.6+0.4,7.0) ; 
    \draw (9.6+2,7.2) -- (9.6+2,7.0) ; 
    \draw (9.6+0.4,7.0) -- (9.6+2,7.0) ; 
    \draw (9.6+1.2,7.0) -- (9.6+1.2,6.8) ; 
    \draw (9.6+1.2-0.2,6.8) -- (9.6+1.2+0.2,6.8) -- (9.6+1.2+0.2,6.4) -- (9.6+1.2-0.2,6.4) -- cycle;
    \draw (9.6+1.2,6.35) node[above] {\footnotesize 0} ;

    \draw (1.2,6.4) -- (1.2,6.2) ; 
    \draw (4.4,6.4) -- (4.4,6.2) ; 
    \draw (1.2,6.2) -- (4.4,6.2) ; 
    \draw (2.8,6.2) -- (2.8,6) ; 
    \draw (2.8-0.2,6) -- (2.8+0.2,6) -- (2.8+0.2,5.6) -- (2.8-0.2,5.6) -- cycle;
    \draw (2.8,5.55) node[above] {\footnotesize 3} ;

    \draw (6.4+1.2,6.4) -- (6.4+1.2,6.2) ; 
    \draw (6.4+4.4,6.4) -- (6.4+4.4,6.2) ; 
    \draw (6.4+1.2,6.2) -- (6.4+4.4,6.2) ; 
    \draw (6.4+2.8,6.2) -- (6.4+2.8,6) ; 
    \draw (6.4+2.8-0.2,6) -- (6.4+2.8+0.2,6) -- (6.4+2.8+0.2,5.6) -- (6.4+2.8-0.2,5.6) -- cycle;
    \draw (6.4+2.8,5.55) node[above] {\footnotesize 0} ;

    \draw (2.8,5.6) -- (2.8,5.4) ; 
    \draw (9.2,5.6) -- (9.2,5.4) ; 
    \draw (2.8,5.4) -- (9.2,5.4) ; 
    \draw (6,5.4) -- (6,5.2) ; 
    \draw (6-0.2,5.2) -- (6+0.2,5.2) -- (6+0.2,4.8) -- (6-0.2,4.8) -- cycle;
    \draw (6,4.75) node[above] {\footnotesize 2} ;

    \draw (-0.2,8.16) node[left] {$\mathfrak{e}_0$} ;
    \draw (0.6,7.36) node[right] {$\mathfrak{e}_1$} ;
    \draw (1.4,6.58) node[right] {$\mathfrak{e}_2$} ;
    \draw (3,5.76) node[right] {$\mathfrak{e}_3$} ;
    \draw (6.2,4.96) node[right] {$\mathfrak{e}_4$} ;

\end{tikzpicture}

 \caption{\leftskip=1.8truecm \rightskip=1.8truecm An example of $\T_4$ with $m=2$. Each vertex is represented by a square box; the number in the box represents the value of the Derrida--Retaux system at the vertex.}

 \label{f:fig1}
\end{figure}

For $v\in \T$ with $|v|=0$ and integer $n\ge 0$, we write $v_n$ for the unique descendant of $v$ in generation $n$, and call $(v=v_0, \, v_1, \, v_2, \ldots, \, v_n)$ the path in $\T$ from $v$ leading to $v_n$;\footnote{Degenerate case: when $|v|=0$, the path from $v$ to $v$ is reduced to the singleton $v$.} for any integer $n\ge 0$, let $\mathtt{bro}(v_n)$ be the set of the ``brothers" of $v_n$, i.e., the parents of $v_{n+1}$ that are not $v_n$. Let
\begin{equation}
    \xi_n(v)
    :=
    \sum_{y\in \mathtt{bro}(v_n)}X(y).
    \label{xi}
\end{equation}

\noindent The path $(v=v_0, \, v_1, \, v_2, \ldots, \, v_n)$ is called {\it open} if 
$$
X(v) +\xi_0(v) + \xi_1(v) + \cdots + \xi_i(v)  \ge i+1, \qquad \forall 0\le i\le n-1.
$$

\noindent Equivalently, the path is open if for $u:= v_i$ (for any $1\le i\le n$), we have $X(u) = X(u^{(1)}) +\cdots+ X(u^{(m)}) -1$ (with $u^{(1)}$, $\ldots$, $u^{(m)}$ denoting as before the parents of $u$),\footnote{Degenerate case: for $|v|=0$, the path consisting of the singleton $v$ is considered as open.} i.e., $X(u^{(1)}) +\cdots+ X(u^{(m)}) \ge 1$.

For $v\in \T$ and integer $i\ge 0$, let 
\begin{eqnarray}
    N^{(i)}(v)
 &:=& \hbox{number of open paths leading to $v$,}
    \nonumber
    \\
 && \hbox{and starting from a vertex $u$ with $|u|=0$ such that $X(u)=i$.}
    \label{N(v)}
\end{eqnarray}

\noindent In particular, $N^{(i)}(v) = {\bf 1}_{\{ X(v)=i\} }$ if $|v|=0$.

For $n\ge 0$, we write
\begin{equation}
    N_n^{(i)} := N^{(i)} (\mathfrak{e}_n),
    \qquad
    N_n := \sum_{i=0}^\infty N_n^{(i)},
    \qquad
    X_n := X(\mathfrak{e}_n) \, .
    \label{N}
\end{equation}

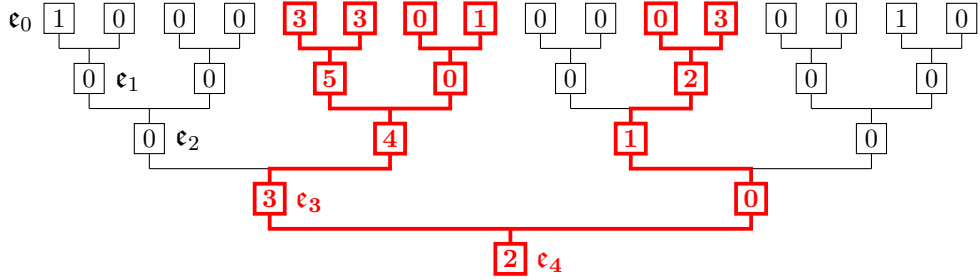
\begin{figure}[ht]
\centering 

\begin{tikzpicture}[scale=1]

    \draw (0,8) -- (0,7.8) ; 
    \draw (-0.2,8) -- (0.2,8) -- (0.2,8.4) -- (-0.2,8.4) -- cycle;
    \draw (0,7.94) node[above] {\footnotesize 1} ;
    \draw (0.8,8) -- (0.8,7.8) ; 
    \draw (0.8-0.2,8) -- (0.8+0.2,8) -- (0.8+0.2,8.4) -- (0.8-0.2,8.4) -- cycle;
    \draw (0.8,7.94) node[above] {\footnotesize 0} ;
    \draw (0,7.8) -- (0.8,7.8) ; 
    \draw (0.4,7.8) -- (0.4,7.6) ; 
    \draw (0.4-0.2,7.6) -- (0.4+0.2,7.6) -- (0.4+0.2,7.2) -- (0.4-0.2,7.2) -- cycle;
    \draw (0.4,7.14) node[above] {\footnotesize 0} ;

    \draw (1.6,8) -- (1.6,7.8) ; 
    \draw (1.6-0.2,8) -- (1.6+0.2,8) -- (1.6+0.2,8.4) -- (1.6-0.2,8.4) -- cycle;
    \draw (1.6,7.94) node[above] {\footnotesize 0} ;
    \draw (1.6+0.8,8) -- (1.6+0.8,7.8) ; 
    \draw (1.6+0.8-0.2,8) -- (1.6+0.8+0.2,8) -- (1.6+0.8+0.2,8.4) -- (1.6+0.8-0.2,8.4) -- cycle;
    \draw (1.6+0.8,7.94) node[above] {\footnotesize 0} ;
    \draw (1.6,7.8) -- (1.6+0.8,7.8) ; 
    \draw (1.6+0.4,7.8) -- (1.6+0.4,7.6) ; 
    \draw (1.6+0.4-0.2,7.6) -- (1.6+0.4+0.2,7.6) -- (1.6+0.4+0.2,7.2) -- (1.6+0.4-0.2,7.2) -- cycle;
    \draw (1.6+0.4,7.14) node[above] {\footnotesize 0} ;

    \draw [line width=0.5mm, red] (3.2,8) -- (3.2,7.8) ; 
    \draw [line width=0.5mm, red] (3.2-0.2,8) -- (3.2+0.2,8) -- (3.2+0.2,8.4) -- (3.2-0.2,8.4) -- cycle;
    \draw (3.2,7.94) node[above] {\color{red}\bf\footnotesize 3} ;
    \draw [line width=0.5mm, red] (3.2+0.8,8) -- (3.2+0.8,7.8) ; 
    \draw [line width=0.5mm, red] (3.2+0.8-0.2,8) -- (3.2+0.8+0.2,8) -- (3.2+0.8+0.2,8.4) -- (3.2+0.8-0.2,8.4) -- cycle;
    \draw (3.2+0.8,7.94) node[above] {\color{red}\bf\footnotesize 3} ;
    \draw [line width=0.5mm, red] (3.2-0.025,7.8) -- (3.2+0.8+0.025,7.8) ; 
    \draw [line width=0.5mm, red] (3.2+0.4,7.8) -- (3.2+0.4,7.6) ; 
    \draw [line width=0.5mm, red] (3.2+0.4-0.2,7.6) -- (3.2+0.4+0.2,7.6) -- (3.2+0.4+0.2,7.2) -- (3.2+0.4-0.2,7.2) -- cycle;
    \draw (3.2+0.4,7.14) node[above] {\color{red}\bf\footnotesize 5} ;

    \draw [line width=0.5mm, red] (4.8,8) -- (4.8,7.8) ; 
    \draw [line width=0.5mm, red] (4.8-0.2,8) -- (4.8+0.2,8) -- (4.8+0.2,8.4) -- (4.8-0.2,8.4) -- cycle;
    \draw (4.8,7.94) node[above] {\color{red}\bf\footnotesize 0} ;
    \draw [line width=0.5mm, red] (4.8+0.8,8) -- (4.8+0.8,7.8) ; 
    \draw [line width=0.5mm, red] (4.8+0.8-0.2,8) -- (4.8+0.8+0.2,8) -- (4.8+0.8+0.2,8.4) -- (4.8+0.8-0.2,8.4) -- cycle;
    \draw (4.8+0.8,7.94) node[above] {\color{red}\bf\footnotesize 1} ;
    \draw [line width=0.5mm, red] (4.8-0.025,7.8) -- (4.8+0.8+0.025,7.8) ; 
    \draw [line width=0.5mm, red] (4.8+0.4,7.8) -- (4.8+0.4,7.6) ; 
    \draw [line width=0.5mm, red] (4.8+0.4-0.2,7.6) -- (4.8+0.4+0.2,7.6) -- (4.8+0.4+0.2,7.2) -- (4.8+0.4-0.2,7.2) -- cycle;
    \draw (4.8+0.4,7.14) node[above] {\color{red}\bf\footnotesize 0} ;

    \draw (6.4,8) -- (6.4,7.8) ; 
    \draw (6.4-0.2,8) -- (6.4+0.2,8) -- (6.4+0.2,8.4) -- (6.4-0.2,8.4) -- cycle;
    \draw (6.4,7.94) node[above] {\footnotesize 0} ;
    \draw (6.4+0.8,8) -- (6.4+0.8,7.8) ; 
    \draw (6.4+0.8-0.2,8) -- (6.4+0.8+0.2,8) -- (6.4+0.8+0.2,8.4) -- (6.4+0.8-0.2,8.4) -- cycle;
    \draw (6.4+0.8,7.94) node[above] {\footnotesize 0} ;
    \draw (6.4,7.8) -- (6.4+0.8,7.8) ; 
    \draw (6.4+0.4,7.8) -- (6.4+0.4,7.6) ; 
    \draw (6.4+0.4-0.2,7.6) -- (6.4+0.4+0.2,7.6) -- (6.4+0.4+0.2,7.2) -- (6.4+0.4-0.2,7.2) -- cycle;
    \draw (6.4+0.4,7.14) node[above] {\footnotesize 0} ;

    \draw [line width=0.5mm, red] (8.0,8) -- (8.0,7.8) ; 
    \draw [line width=0.5mm, red] (8.0-0.2,8) -- (8.0+0.2,8) -- (8.0+0.2,8.4) -- (8.0-0.2,8.4) -- cycle;
    \draw (8.0,7.94) node[above] {\color{red}\bf\footnotesize 0} ;
    \draw [line width=0.5mm, red] (8.0+0.8,8) -- (8.0+0.8,7.8) ; 
    \draw [line width=0.5mm, red] (8.0+0.8-0.2,8) -- (8.0+0.8+0.2,8) -- (8.0+0.8+0.2,8.4) -- (8.0+0.8-0.2,8.4) -- cycle;
    \draw (8.0+0.8,7.94) node[above] {\color{red}\bf\footnotesize 3} ;
    \draw [line width=0.5mm, red] (8.0-0.025,7.8) -- (8.0+0.8+0.025,7.8) ; 
    \draw [line width=0.5mm, red] (8.0+0.4,7.8) -- (8.0+0.4,7.6) ; 
    \draw [line width=0.5mm, red] (8.0+0.4-0.2,7.6) -- (8.0+0.4+0.2,7.6) -- (8.0+0.4+0.2,7.2) -- (8.0+0.4-0.2,7.2) -- cycle;
    \draw (8.0+0.4,7.14) node[above] {\color{red}\bf\footnotesize 2} ;

    \draw (9.6,8) -- (9.6,7.8) ; 
    \draw (9.6-0.2,8) -- (9.6+0.2,8) -- (9.6+0.2,8.4) -- (9.6-0.2,8.4) -- cycle;
    \draw (9.6,7.94) node[above] {\footnotesize 0} ;
    \draw (9.6+0.8,8) -- (9.6+0.8,7.8) ; 
    \draw (9.6+0.8-0.2,8) -- (9.6+0.8+0.2,8) -- (9.6+0.8+0.2,8.4) -- (9.6+0.8-0.2,8.4) -- cycle;
    \draw (9.6+0.8,7.94) node[above] {\footnotesize 0} ;
    \draw (9.6,7.8) -- (9.6+0.8,7.8) ; 
    \draw (9.6+0.4,7.8) -- (9.6+0.4,7.6) ; 
    \draw (9.6+0.4-0.2,7.6) -- (9.6+0.4+0.2,7.6) -- (9.6+0.4+0.2,7.2) -- (9.6+0.4-0.2,7.2) -- cycle;
    \draw (9.6+0.4,7.14) node[above] {\footnotesize 0} ;

    \draw (11.2,8) -- (11.2,7.8) ; 
    \draw (11.2-0.2,8) -- (11.2+0.2,8) -- (11.2+0.2,8.4) -- (11.2-0.2,8.4) -- cycle;
    \draw (11.2,7.94) node[above] {\footnotesize 1} ;
    \draw (11.2+0.8,8) -- (11.2+0.8,7.8) ; 
    \draw (11.2+0.8-0.2,8) -- (11.2+0.8+0.2,8) -- (11.2+0.8+0.2,8.4) -- (11.2+0.8-0.2,8.4) -- cycle;
    \draw (11.2+0.8,7.94) node[above] {\footnotesize 0} ;
    \draw (11.2,7.8) -- (11.2+0.8,7.8) ; 
    \draw (11.2+0.4,7.8) -- (11.2+0.4,7.6) ; 
    \draw (11.2+0.4-0.2,7.6) -- (11.2+0.4+0.2,7.6) -- (11.2+0.4+0.2,7.2) -- (11.2+0.4-0.2,7.2) -- cycle;
    \draw (11.2+0.4,7.14) node[above] {\footnotesize 0} ;

    \draw (0.4,7.2) -- (0.4,7.0) ; 
    \draw (2,7.2) -- (2,7.0) ; 
    \draw (0.4,7.0) -- (2,7.0) ; 
    \draw (1.2,7.0) -- (1.2,6.8) ; 
    \draw (1.2-0.2,6.8) -- (1.2+0.2,6.8) -- (1.2+0.2,6.4) -- (1.2-0.2,6.4) -- cycle;
    \draw (1.2,6.35) node[above] {\footnotesize 0} ;

    \draw [line width=0.5mm, red] (3.2+0.4,7.2) -- (3.2+0.4,7.0) ; 
    \draw [line width=0.5mm, red] (3.2+2,7.2) -- (3.2+2,7.0) ; 
    \draw [line width=0.5mm, red] (3.2+0.4-0.025,7.0) -- (3.2+2+0.025,7.0) ; 
    \draw [line width=0.5mm, red] (3.2+1.2,7.0) -- (3.2+1.2,6.8) ; 
    \draw [line width=0.5mm, red] (3.2+1.2-0.2,6.8) -- (3.2+1.2+0.2,6.8) -- (3.2+1.2+0.2,6.4) -- (3.2+1.2-0.2,6.4) -- cycle;
    \draw (3.2+1.2,6.35) node[above] {\color{red}\bf\footnotesize 4} ;

    \draw (6.4+0.4,7.2) -- (6.4+0.4,7.0) ; 
    \draw [line width=0.5mm, red] (6.4+2,7.2) -- (6.4+2,7.0) ; 
    \draw (6.4+0.4,7.0) -- (6.4+1.2,7.0) ; 
    \draw [line width=0.5mm, red] (6.4+1.2-0.025,7.0) -- (6.4+2+0.025,7.0) ; 
    \draw [line width=0.5mm, red] (6.4+1.2,7.0) -- (6.4+1.2,6.8) ; 
    \draw [line width=0.5mm, red] (6.4+1.2-0.2,6.8) -- (6.4+1.2+0.2,6.8) -- (6.4+1.2+0.2,6.4) -- (6.4+1.2-0.2,6.4) -- cycle;
    \draw (6.4+1.2,6.35) node[above] {\color{red}\bf\footnotesize 1} ;

    \draw (9.6+0.4,7.2) -- (9.6+0.4,7.0) ; 
    \draw (9.6+2,7.2) -- (9.6+2,7.0) ; 
    \draw (9.6+0.4,7.0) -- (9.6+2,7.0) ; 
    \draw (9.6+1.2,7.0) -- (9.6+1.2,6.8) ; 
    \draw (9.6+1.2-0.2,6.8) -- (9.6+1.2+0.2,6.8) -- (9.6+1.2+0.2,6.4) -- (9.6+1.2-0.2,6.4) -- cycle;
    \draw (9.6+1.2,6.35) node[above] {\footnotesize 0} ;

    \draw (1.2,6.4) -- (1.2,6.2) ; 
    \draw [line width=0.5mm, red] (4.4,6.4) -- (4.4,6.2) ; 
    \draw (1.2,6.2) -- (2.8,6.2) ; 
    \draw [line width=0.5mm, red] (2.8-0.025,6.2) -- (4.4+0.025,6.2) ; 
    \draw [line width=0.5mm, red] (2.8,6.2) -- (2.8,6) ; 
    \draw [line width=0.5mm, red] (2.8-0.2,6) -- (2.8+0.2,6) -- (2.8+0.2,5.6) -- (2.8-0.2,5.6) -- cycle;
    \draw (2.8,5.55) node[above] {\color{red}\bf\footnotesize 3} ;

    \draw [line width=0.5mm, red] (6.4+1.2,6.4) -- (6.4+1.2,6.2) ; 
    \draw (6.4+4.4,6.4) -- (6.4+4.4,6.2) ; 
    \draw [line width=0.5mm, red] (6.4+1.2-0.025,6.2) -- (6.4+2.8+0.025,6.2) ; 
    \draw (6.4+2.8,6.2) -- (6.4+4.4,6.2) ; 
    \draw [line width=0.5mm, red] (6.4+2.8,6.2) -- (6.4+2.8,6) ; 
    \draw [line width=0.5mm, red] (6.4+2.8-0.2,6) -- (6.4+2.8+0.2,6) -- (6.4+2.8+0.2,5.6) -- (6.4+2.8-0.2,5.6) -- cycle;
    \draw (6.4+2.8,5.55) node[above] {\color{red}\bf\footnotesize 0} ;

    \draw [line width=0.5mm, red] (2.8,5.6) -- (2.8,5.4) ; 
    \draw [line width=0.5mm, red] (9.2,5.6) -- (9.2,5.4) ; 
    \draw [line width=0.5mm, red] (2.8-0.025,5.4) -- (9.2+0.025,5.4) ; 
    \draw [line width=0.5mm, red] (6,5.4) -- (6,5.2) ; 
    \draw [line width=0.5mm, red] (6-0.2,5.2) -- (6+0.2,5.2) -- (6+0.2,4.8) -- (6-0.2,4.8) -- cycle;
    \draw (6,4.75) node[above] {\color{red}\bf\footnotesize 2} ;

    \draw (-0.2,8.16) node[left] {$\mathfrak{e}_0$} ;
    \draw (0.6,7.36) node[right] {$\mathfrak{e}_1$} ;
    \draw (1.4,6.58) node[right] {$\mathfrak{e}_2$} ;
    \draw (3,5.76) node[right] {$\color{red}\bf\mathfrak{e}_3$} ;
    \draw (6.2,4.96) node[right] {$\color{red}\bf\mathfrak{e}_4$} ;

\end{tikzpicture}

 \caption{\leftskip=1.6truecm \rightskip=1.6truecm The same example as in Figure \ref{f:fig1}. Open paths from the initial generation to $\mathfrak{e}_4$ are marked in bold (and coloured in red), with $N_4^{(0)} = 2$, $N_4^{(1)} =1$, $N_4^{(2)} =0$, $N_4^{(3)}=3$, $N_4^{(i)}=0$ if $i\ge 4$, and $N_4=6$.}

 \label{f:fig2}
\end{figure}

\noindent [The identification of $X_n$ with $X(\mathfrak{e}_n)$ should cause no confusion because only the individual law of each $X_n$ matters.] See Figure \ref{f:fig2} for an example.

\bigskip

\medskip

\begin{definition}
\label{def:pivot}

Let $v$ be a vertex with $|v|=0$. For $k$, $n\ge 0$ and $\ell\ge 0$, we say that $v$ is a $(k, \, \ell)_n$--{\bf pivotal vertex} if $k+\xi_0(v) + \xi_1(v) +\cdots+\xi_{n-1}(v)=n+\ell$ and
$$
k+\xi_0(v) + \xi_1(v) + \cdots + \xi_i(v)  \ge i+1, \qquad \forall 0\le i\le n-1.
$$

 
\end{definition}

\medskip

By definition, for a vertex $v$ in the initial generation, whether it is pivotal or not does not depend on the value of $X(v)$.

The notion of pivotal vertices and that of open paths are obviously related. Let $v$ be a vertex in the initial generation. If $v$ is $(k, \, \ell)_n$--pivotal and $X(v)=k$, then the path from $v$ to $v_n$ (the descendant of $v$ at generation $n$) is open, and $X(v_n)=\ell$. Conversely, if the path from $v$ to $v_n$ is open, then $v$ is $(k, \, \ell)_n$--pivotal with $k:= X(v)$ and $\ell := X(v_n)$. 

Define
\begin{equation}
    A_n^{(k,\ell)}
   :=
    \{ v\in \mathbb{T}_n: \; |v|=0, \; v \hbox{ \rm is $(k, \, \ell)_n$--pivotal} \} ,
    \label{pivot}
\end{equation}

\noindent which stands for the set of $(k, \, \ell)_n$--pivotal vertices in the initial generation of $\mathbb{T}_n$. Then
\begin{equation}
    A_n^{(k,\ell_1)} \cap A_n^{(k,\ell_2)}
    =
    \varnothing \quad \hbox{\rm for all $k\ge 0$ and $\ell_1>\ell_2 \ge 0$.}
    \label{e:Akldisjoint}
\end{equation}

\noindent Define, for integers $n\ge 0$, $k\ge 0$ and $\ell \ge 0$,
\begin{eqnarray}
    \mu_n(k)
 &:=& m^{-k}\, \E([1-(m-1)X_n]m^{X_n} \, {\bf 1}_{\{X_n\le k\}}) \, ,
    \label{mu}
    \\
    S_n^{(k,\ell)} 
 &:=& |A_n^{(k,\ell)}| \, ,
    \label{S}
\end{eqnarray}

\noindent with $|A_n^{(k,\ell)}|$ denoting the cardinality of $A_n^{(k,\ell)}$. We have observed that the notions of pivotal vertices and open paths are related; let us give a more quantitative relation.
Since $\mu_0(0) = \P(X_0=0)$, we have, for $n\ge 0$ and $\ell \ge 0$,
\begin{equation}
    \mu_0(0) \, \E(S_n^{(0,\ell)}) 
    =
    \E(N_n^{(0)} \, {\bf 1}_{\{ X_n =\ell\} }) \, .
    \label{rel_pivot_openpaths}
\end{equation}

Our main result in this section is the following identity. 

\medskip

\begin{theorem}
 \label{t:new_iteration} 
 
 For $n\ge 0$ and $\ell \ge 0$,
 \begin{equation}
     \sum_{k=0}^\infty \mu_0(k) \, \E(S_n^{(k,\ell)})
     =
     \mu_n(\ell) \, .
     \label{new_iteration}
 \end{equation}

\end{theorem}

\medskip

\noindent {\it Proof.} We prove \eqref{new_iteration} by induction on $n$. When $n=0$, \eqref{new_iteration} holds since $S_0^{(k,\ell)} = {\bf 1}_{\{k=\ell\}}$. Assume \eqref{new_iteration} holds for $n$. Let us prove it for $n+1$. By the definition of pivotal vertices, for any $k$, $\ell\ge 0$ and $v\in \T_n$ with $|v|=0$,
\begin{equation}
    {\bf 1}_{\{v\in A_{n+1}^{(k,\ell)}\}}
    =
    \sum_{j=0}^{\ell+1} {\bf 1}_{\{\xi_n(v)=j\}} \, {\bf 1}_{\{v\in A_n^{(k, \ell+1-j)}\}} \, .
    \label{A(n+1)_and_An}
\end{equation}

\noindent By symmetry, $\E(S_{n+1}^{(k,\ell)}) = \sum_{v\in \T_{n+1}: |v|=0} \E({\bf 1}_{\{v\in A_{n+1}^{(k,\ell)}\}} ) =  m^{n+1} \E( {\bf 1}_{\{\mathfrak{e}_0\in A_{n+1}^{(k,\ell)}\}})$. Putting these identities together, we get
$$
\E(S_{n+1}^{(k,\ell)})
=
m^{n+1} \sum_{j=0}^{\ell+1}\E( {\bf 1}_{\{\xi_n(\mathfrak{e}_0)=j\}} \, {\bf 1}_{\{\mathfrak{e}_0\in A_n^{(k, \ell+1-j)}\}}) \, .
$$

\noindent Since $\{\xi_n(\mathfrak{e}_0)=j\}$ and $\{\mathfrak{e}_0\in A_n^{(k, \ell+1-j)}\}$ are independent events, this leads to:
\begin{eqnarray*}
    \E(S_{n+1}^{(k,\ell)})
 &=& m^{n+1} \sum_{j=0}^{\ell+1} \P(\xi_n(\mathfrak{e}_0)=j) \, \E({\bf 1}_{\{\mathfrak{e}_0\in A_n^{(k,\ell+1-j)}\}})
    \\
 &=& m\sum_{j=0}^{\ell+1}\P(\xi_n(\mathfrak{e}_0)=j) \E(S_n^{(k,\ell+1-j)}) \, ,
\end{eqnarray*}

\noindent which implies that\footnote{We can exchange the order of the sums because $\sum_{k=0}^\infty$ is a finite sum: $\E(S_{n+1}^{(k,\ell)}) =0$ for $k>(n+1)+\ell$.}
\begin{eqnarray*}
    \sum_{k=0}^\infty \mu_0(k) \E(S_{n+1}^{(k,\ell)})
 &=& m\sum_{k=0}^\infty \mu_0(k) \sum_{j=0}^{\ell+1} \P(\xi_n(\mathfrak{e}_0)=j) \, \E(S_n^{(k,\ell+1-j)})
    \\
 &=& m\sum_{j=0}^{\ell+1} \P(\xi_n(\mathfrak{e}_0)=j)\sum_{k=0}^\infty \mu_0(k) \, \E(S_n^{(k,\ell+1-j)}).
\end{eqnarray*}

\noindent By the induction assumption, we get that $\sum_{k=0}^\infty \mu_0(k) \E(S_n^{(k,\ell+1-j)}) = \mu_n(\ell+1-j)$ (for $j\le \ell+1$), which is $m^{-(\ell+1-j)} \, \E [(1-(m-1)X_n)m^{X_n}\, {\bf 1}_{\{ X_n\le \ell+1-j\} }]$. So we have
\begin{eqnarray}
    \sum_{k=0}^\infty \mu_0(k) \E(S_{n+1}^{(k,\ell)})
 &=& m^{-\ell}\sum_{j=0}^{\ell+1} m^j \,\P(\xi_n(\mathfrak{e}_0)=j) \, \E [(1-(m-1)X_n)m^{X_n}\, {\bf 1}_{\{ X_n\le \ell+1-j\} }]
    \nonumber 
    \\
 &=& m^{-\ell} \, \E [(1-(m-1)X_n) m^{X_n+\xi_n(\mathfrak{e}_0)}\, {\bf 1}_{\{ X_n+\xi_n(\mathfrak{e}_0) \le \ell+1\} }]
    \nonumber
    \\
 &=& m^{-\ell} \, \E [(1-(m-1)X_n^{(1)}) m^{X_n^{(1)}+\cdots+X_n^{(m)}}\, {\bf 1}_{\{ X_n^{(1)}+\cdots+X_n^{(m)}\le \ell+1\} }] \, ,
    \label{e:NYxiequal}
\end{eqnarray}

\noindent the last equality being a consequence of the fact that the two-dimensional random vector $(X_n, \, \xi_n(\mathfrak{e}_0))$, which is $(X(\mathfrak{e}_n), \, \xi_n(\mathfrak{e}_0))$ by definition, has the same distribution as $(X_n^{(1)}, \, X_n^{(2)}+\cdots+X_n^{(m)})$. 

On the other hand, by symmetry,
\begin{eqnarray*}
 && \E [(1-(m-1)X_{n+1})m^{X_{n+1}}\, {\bf 1}_{\{ X_{n+1}\le \ell\} }]
    -
    \P(X_n^{(1)}+\cdots+X_n^{(m)}=0)
    \\
 &=& \sum_{u=1}^{\ell+1} \E [(1-(m-1)(X_n^{(1)}+\cdots+X_n^{(m)} -1)) m^{u-1}\, {\bf 1}_{\{ X_n^{(1)}+\cdots+X_n^{(m)}=u\} }]
    \\
 &=& \sum_{u=1}^{\ell+1} \E [(1-(m-1)(m X_n^{(1)} -1)) m^{u-1}\, {\bf 1}_{\{ X_n^{(1)}+\cdots+X_n^{(m)}=u\} }]
    \\
 &=& \sum_{u=1}^{\ell+1} \E [(m-m(m-1) X_n^{(1)}) m^{u-1}\, {\bf 1}_{\{ X_n^{(1)}+\cdots+X_n^{(m)}=u\} }]
    \\
 &=& \sum_{u=1}^{\ell+1} \E [(1-(m-1)X_n^{(1)}) m^u\, {\bf 1}_{\{ X_n^{(1)}+\cdots+X_n^{(m)}=u\} }] \, .
\end{eqnarray*}

\noindent Hence
\begin{eqnarray*}
 &&\E [(1-(m-1)X_{n+1})m^{X_{n+1}}\, {\bf 1}_{\{ X_{n+1}\le \ell\} }]
    \\
 &=& \E [(1-(m-1)X_n^{(1)}) m^{X_n^{(1)}+\cdots+X_n^{(m)}}\, {\bf 1}_{\{ X_n^{(1)}+\cdots+X_n^{(m)}\le \ell+1\} }] \, .
\end{eqnarray*}

\noindent Going back to \eqref{e:NYxiequal}, we get
$$
\sum_{k=0}^\infty \mu_0(k)\, \E(S_{n+1}^{(k,\ell)})
=
m^{-\ell}\, \E [(1-(m-1)X_{n+1})m^{X_{n+1}}\, {\bf 1}_{\{ X_{n+1}\le \ell\} }],
$$

\noindent which is the desired inequality \eqref{new_iteration} for $n+1$. By induction, \eqref{new_iteration} holds for all $n\ge 0$. Theorem \ref{t:new_iteration} is proved.\qed

\medskip

\begin{remark}

 When the Derrida--Retaux system is critical or subcritical, Theorem \ref{t:new_iteration} has a simple consequence that is useful in the proof of the upper bound in Theorem \ref{t:main}. [Its most important application, however, will arrive in Lemma \ref{l:E(S)>} as a technical ingredient in the proof of the lower bound in Theorem \ref{t:main}.] Assume $(X_n, \, n\ge 0)$ is critical or subcritical: $(m-1) \E(X_0 \, m^{X_0}) \le \E(m^{X_0})<\infty$. This yields (see \cite{collet-eckmann-glaser-martin}, \cite{bmxyz_questions}; a simple explanation is given at the beginning of Section \ref{s:preparation_ub}) that $(m-1) \E(X_n \, m^{X_n}) \le \E(m^{X_n})<\infty$ for all $n\ge 0$. In particular, 
 $$
 \mu_n(k)
 \ge
 m^{-k}\, \E([(m-1)X_n-1]m^{X_n} \, {\bf 1}_{\{X_n\ge k+1\}})
 \ge 0,
 $$
 
 \noindent for all $k\ge 0$. On the other hand, by definition, 
 $$
 \mu_n(k) 
 \le 
 m^{-k}\, \E([1-(m-1)X_n]m^{X_n} \, {\bf 1}_{\{X_n=0\}}) 
 = 
 m^{-k} \, \P(X_n=0) 
 \le 
 m^{-k} \, .
 $$
 
 \noindent So for critical or subcritical Derrida--Retaux systems, we have
 \begin{equation}
     0\le \mu_n(k) \le m^{-k} ,
     \qquad
     k\ge 0, \; n\ge 0\, .
     \label{mu_bounds}
 \end{equation}

 Let us look at Theorem \ref{t:new_iteration}. By considering only the term $k=0$, and since $\mu_n(\, \cdot\,)$ is non-negative, the theorem implies that for integers $n\ge 0$ and $\ell \ge 0$,
 $$
 \mu_0(0) \, \E(S_n^{(0,\ell)})
 \le
 \mu_n(\ell)
 \le 
 m^{-\ell} \, ,
 $$

 \noindent the last inequality being from \eqref{mu_bounds}. Let $N_n^{(0)}$ be as in \eqref{N}. We have $\E(N_n^{(0)} \, {\bf 1}_{\{ X_n =\ell\} }) = \mu_0(0) \, \E(S_n^{(0,\ell)})$ (see \eqref{rel_pivot_openpaths}). Hence
 $$
 \E(N_n^{(0)}\, {\bf 1}_{\{ X_n = \ell\} })
 \le 
 m^{-\ell} \, .
 $$

 \noindent Consequently,
 \begin{equation}
     \E(N_n^{(0)}\, {\bf 1}_{\{ X_n \ge \ell\} })
     \le
     \frac{c_1}{m^\ell} ,
     \qquad
     n\ge 0, \; \ell\ge 0 \, ,
     \label{E(N_Y>k)}
 \end{equation}
 
 \noindent with $c_1 := \frac{m}{m-1}$.\qed
   
\end{remark}

\begin{remark}

 Let $N_n$ denote, as in \eqref{N}, the number of open paths from the initial generation to $\mathfrak{e}_n$. When $N_n \ge 1$, we say that the system is sustainable until $\mathfrak{e}_n$. We are going to see that in the study of the Derrida--Retaux system, it is important to understand the situation of sustainability (more important than to understand the situation $X_n \ge 1$).
 
 By definition,
 $$
 \{ X_n \ge 1\} \subset \{ N_n \ge 1\} = \{ X_{n-1}^{(1)} + \cdots + X_{n-1}^{(m)} \ge 1 \} ,
 $$
 
 \noindent where $X_{n-1}^{(1)}$, $\ldots$, $X_{n-1}^{(m)}$ are, as before, independent copies of $X_{n-1}$. Therefore,
 $$
 \P(X_n \ge 1) \le \P (N_n \ge 1) \le m\, \P(X_{n-1} \ge 1) \, .
 $$
 
 \noindent In particular, for a critical system $(Y_n, \, n\ge 0)$, if $\E(t^{Y_0})<\infty$ (for some $t>m$), then Theorem \ref{t:main} says that the probability that the system is sustainable until $\mathfrak{e}_n$ equals $\frac{1}{n^{2+o(1)}}$ (for $n\to \infty$).\qed

\end{remark}

\section{Preliminaries on the moment generating function}
\label{s:preparation_ub}

This section collects some estimates, in four lemmas below, for the moment generating function of Derrida--Retaux systems at criticality. All the results are either known, or are proved painlessly by means of existing methods, with the only exception of Lemma \ref{l:YLNUPP}~(i) which is obtained as a consequence of \eqref{E(N_Y>k)}.   

Let $(X_n, \, n\ge 0)$ be a Derrida--Retaux system defined via \eqref{iteration} (it will soon be assumed to be critical), satisfying $\E(X_0 \, m^{X_0})<\infty$. Write
$$
G_n(s) := \E(s^{X_n}),
$$

\noindent for the moment generating function of $X_n$. The recursion in \eqref{iteration} means
\begin{equation}
    G_{n+1}(s) 
    =
    \frac1s \, G_n(s)^m + \Big( 1- \frac1s\Big) \, G_n(0)^m ,
    \label{iteration_Gn}
\end{equation}

\noindent from which we get 
$$
G'_{n+1}(s) 
=
\frac{m}{s} \, G_n'(s) G_n(s)^{m-1} - \frac{1}{s^2} \, G_n(s)^m + \frac{1}{s^2} \, G_n(0)^m \, .
$$

\noindent Hence 
\begin{equation}
    s(s-1) G'_{n+1}(s) - G_{n+1}(s)
    =
    [m(s-1) G'_n(s) - G_n(s)] G_n(s)^{m-1} \, .
    \label{p11}
\end{equation}

\noindent In particular,
$$
m(m-1) G'_{n+1}(m) - G_{n+1}(m)
=
[m(m-1) G'_n(m) - G_n(m)] G_n(m)^{m-1} \, .
$$

\noindent Iterating the identity leads to the following formula: for $n\ge 1$,
$$
m(m-1) G'_n(m) - G_n(m)
=
[m(m-1) G'_0(m) - G_0(m)] \prod_{j=0}^{n-1} G_j(m)^{m-1} \, .
$$

\noindent In particular, it tells us that the sign of $m(m-1) G'_n(m) - G_n(m)$ stays the same for all $n\ge 0$. For example, if $(m-1)\E(X_0 \, m^{X_0}) \le \E(m^{X_0})$, then for all $n\ge 0$, $(m-1)\E(X_n \, m^{X_n}) \le \E(m^{X_n})$; and since $\E(X_n \, m^{X_n}) \ge \E(X_n) \, \E(m^{X_n})$ (by the FKG inequality), we get $\E(X_n) \le \frac{1}{m-1}$, $\forall n\ge 0$, which is the second part of Theorem A recalled in the introduction. Another consequence is that if $(X_n, \, n\ge 0)$ is critical (resp.~supercritical; subcritical), then so is $(X_{n+k}, \, n\ge 0)$ for any integer $k\ge 0$.

We can extend \eqref{iteration_Gn} by bringing in also the number of open paths. For any vertex $v\in \T$ with $|v| \ge 1$, let $v^{(1)}$, $\ldots$, $v^{(m)}$ denote as before its parents in generation $|v|-1$. By definition of $N^{(i)}(v)$ in \eqref{N(v)}, for any integer $i\ge 0$,
\begin{eqnarray*}
    X(v)
 &=& [ X(v^{(1)}) +\cdots+ X(v^{(m)})-1 ]^+,
    \\
    N^{(i)}(v)
 &=& [ N^{(i)}(v^{(1)}) +\cdots+ N^{(i)}(v^{(m)}) ]\, {\bf 1}_{\{ X(v^{(1)}) +\cdots+ X(v^{(m)}) \ge 1\} } \, .
\end{eqnarray*}

\noindent With $\Sigma_X = \Sigma_X(v) := X(v^{(1)}) + \cdots + X(v^{(m)})$ and $\Sigma_N = \Sigma_N(v, \, i) := N^{(i)}(v^{(1)}) + \cdots + N^{(i)}(v^{(m)})$, this can be rewritten as
$$
(X(v), \, N^{(i)}(v)) 
=
(\Sigma_X -1, \, \Sigma_N) \, {\bf 1}_{ \{ \Sigma_X \ge 1\} } 
+
(0, \, 0) \, {\bf 1}_{ \{ \Sigma_X =0\} } \, .
$$

\noindent Thus for $s> 0$ and $t\ge 0$ (with $0^0 := 1$),
\begin{eqnarray*}
    \E(s^{X(v)} \, t^{N^{(i)}(v)})
 &=& \E(s^{\Sigma_X -1} \, t^{\Sigma_N}\, {\bf 1}_{\{ \Sigma_X \ge 1 \} })
    +
    \P(\Sigma_X =0)
    \\
 &=& \E(s^{\Sigma_X -1} \, t^{\Sigma_N})
    -
    \frac1s \, \E(t^{\Sigma_N}\, {\bf 1}_{\{ \Sigma_X =0 \} })
    +
    \P(\Sigma_X=0) .
\end{eqnarray*}

Define, for any integer $i\ge 0$,
$$
G_n(s, \, t)
=
G_n(s, \, t, \, i)
:=
\E(s^{X_n} \, t^{N^{(i)}_n}) ,
\qquad
s\ge 0, \; t \ge 0,
$$

\noindent which is the joint moment generating function for the pair $(X_n, \, N^{(i)}_n)$ ($N^{(i)}_n$ being from \eqref{N}). [In particular, $G_n(s) = G_n(s, \, 1)$.] Since $(X(v^{(j)}), \, N^{(i)}(v^{(j)}))$, for $1\le j\le m$, are i.i.d.\ having the law of $(X_{n-1}, \, N^{(i)}_{n-1})$, we have 
\begin{eqnarray*}
    \E(s^{\Sigma_X -1} \, t^{\Sigma_N}) 
 &=& \frac1s \, G_{n-1} (s, \, t)^m ,
    \\
    \E(t^{\Sigma_N}\, {\bf 1}_{\{ \Sigma_X =0 \} }) 
 &=& [\E(t^{N_{n-1}}\, {\bf 1}_{\{ X_{n-1} =0 \} })]^m 
    = G_{n-1}(0, \, t)^m,
    \\
    \P(\Sigma_X =0) 
 &=& [\P( X_{n-1} =0)]^m 
    =
    G_{n-1}(0, \, 1)^m\, .
\end{eqnarray*}

\noindent Consequently, for $n\ge 1$, $s> 0$ and $t\ge 0$,
\begin{equation}
    G_n (s, \, t)
    =
    \frac{G_{n-1}(s, \, t)^m}{s}
    -
    \frac{G_{n-1}(0, \, t)^m}{s}
    +
    G_{n-1}(0, \, 1)^m \, .
    \label{pf:recursive_formula3}
\end{equation}

\noindent This formula is useful in the proof of Lemma \ref{l:YLNUPP} below.

In the rest of this section, we deal with the Derrida--Retaux system at criticality, and denote it by $(Y_n, \, n\ge 0)$: $(m-1)\E(Y_n \, m^{Y_n}) = \E(m^{Y_n})$ for all $n\ge 0$. We assume $\E(Y_0^3 \, m^{Y_0})<\infty$, which implies $\E(Y_n^3 \, m^{Y_n})<\infty$ for all $n\ge 0$.


\medskip

\begin{lemma}
 \label{l:G_nk123} 

 Assume $\E(Y_0^3 \, m^{Y_0})<\infty$. Let $G_n(s) := \E(s^{Y_n})$. There exist positive constants $c_2$, $c_3$, $c_4$, $c_5$, $c_6$, $c_7$ and $c_8$ such that for $n\ge 1$,

{\rm (i)} $\E(m^{Y_n})\le c_2$;

{\rm (ii)} $\frac{1}{m-1} \le \E(Y_n \, m^{Y_n})\le c_3$;

{\rm (iii)} $\E(Y_n^2 \, m^{Y_n})\le c_4\, n$;

{\rm (iv)} $\frac{c_5}{n} \le \E(m^{Y_n} \, {\bf 1}_{\{Y_n\ge 1\}}) \le \frac{c_6}{n}$;

{\rm (v)} $c_7 \, n^2\le \prod_{i=0}^{n-1}G_i(m)^{m-1}\le c_8 \, n^2$.

\end{lemma}

\medskip

Lemma \ref{l:G_nk123} is known: (v) was proved in \cite{collet-eckmann-glaser-martin} (case $m=2$) and \cite{bmxyz_questions}, whereas (i) follows from Theorem A~(ii) (it is also a consequence of (v)); since the system is critical, we have $\E(Y_n \, m^{Y_n}) = \frac{1}{m-1} \, \E(m^{Y_n})$, so the second inequality in (ii) follows from (i); for (iii) and (iv), see \cite{bmvxyz_conjecture_DR} and \cite{bmxyz_questions}, respectively.  

The first four parts of Lemma \ref{l:G_nk123} can be stated in a more compact way and may be extended as follows: Suppose $\E(Y_0^j m^{Y_0})<\infty$ for some integer $j \ge 2$. For any $k\in \{0, \, 1, \, ... \, j\}$,  there exist some constants $c_9>0$ and $c_{10}>0$ depending on $k$ \footnote{In Lemma \ref{l:moment_derivative_lb} below, we are going to work under a stronger integrability assumption, and give dependence in $k$ of the constants $c_9$ and $c_{10}$.} such that for all integers $n\ge 1$, 
\begin{equation}
    c_9 \, n^{k-1} \le \E(Y_n^k m^{Y_n} \, {\bf 1}_{\{Y_n\ge 1\}}) \le c_{10} \, n^{k-1} \, .
    \label{Yn_moments}
\end{equation}

\noindent Plainly Lemma \ref{l:G_nk123} corresponds to the case $j=3$. 


Let us prove \eqref{Yn_moments}; the proof of the forthcoming Lemmas \ref{l:YLNUPP} and \ref{l:moment_derivative_lb} will be a refinement of the same argument. Only the case $k\in \{2, \, 3, \, \ldots, \, j\}$ needs to be proved because for $k\in \{0, 1\}$, both inequalities in \eqref{Yn_moments} are already contained in Lemma \ref{l:G_nk123}. By \eqref{iteration_Gn},
$$
s G_{n+1}(s) = G_n(s)^m + (s-1) G_n(0)^m\, .
$$

\noindent On both sides, we differentiate $k$ times with respect to $s$ (since $k \ge 2$, the second term on the right-hand side, being affine in $s$, disappears after the second differentiation), and apply the general Leibniz rule to the first term on the right-hand side, to see that
\begin{equation}
    s G_{n+1}^{(k)}(s) + k G_{n+1}^{(k-1)}(s)
    =
    \sum_{(k_1, \ldots, k_m) \in \z_+^m: \, k_1 + \cdots + k_m = k} \frac{k!}{k_1! \cdots k_m!} \prod_{i=1}^m G_n^{(k_i)}(s) \, .
    \label{recurrence_derivatives_Gn}
\end{equation}

\noindent 
 
 We prove the upper bound in \eqref{Yn_moments} (for all $k \in \{2, 3,  ..., j\}$) by induction in $k$. Note that the case  $k=2$ follows from Lemma \ref{l:G_nk123}  (iii), whereas the case $k=3$ is  known from \cite{bmxyz_questions} (equation (19) and Proposition 1 there). Consider   $k\ge 4$ and  suppose that the upper bound  in \eqref{Yn_moments} holds for all cases up to $k-1$. Taking $s=m$ in   \eqref{recurrence_derivatives_Gn} we get that   $G_{n+1}^{(k)}(m)\le G_n(m)^{m-1}\, G_n^{(k)}(m) + c_{11} n^{k-2} $ for some constant $c_{11}>0$ only depending on $k$.  Iterating the inequality gives that for $n\ge 2$,
$$
G_n^{(k)}(m) 
\le
\Big( G_1^{(k)}(m) + \sum_{\ell=1}^{n-1} \frac{c_{11} \, \ell^{k-2}}{\prod_{i=1}^\ell G_i(m)^{m-1}} \Big) \prod_{i=1}^{n-1} G_i(m)^{m-1} \, .
$$
 
\noindent Using Lemma \ref{l:G_nk123}~(v), this yields the upper bound  in \eqref{Yn_moments}.
 
For the lower bound, we use the Cauchy--Schwarz inequality to see that for $k\ge 2$,
\begin{equation}
    \E(Y_n^k m^{Y_n} \, {\bf 1}_{\{Y_n\ge 1\}}) \ge \frac{[\E(Y_n^{k-1} m^{Y_n} \, {\bf 1}_{\{Y_n\ge 1\}})]^2}{\E(Y_n^{k-2} m^{Y_n} \, {\bf 1}_{\{Y_n\ge 1\}})}. 
    \label{CS}
\end{equation}
 
\noindent This, together with the upper bound in \eqref{Yn_moments} which we have just proved, yield the lower bound by induction in $k$. 



We now proceed to our second lemma, which focuses on the joint moment generating function of $Y_n$ and $N_n^{(0)}$. The lemma will be useful in the proof of the upper bound in Theorem \ref{t:main}.

\medskip

\begin{lemma}
\label{l:YLNUPP}

 Assume $\E(Y_0^3 m^{Y_0})<\infty$. Let $N_n^{(0)}$ denote the number of open paths as in \eqref{N} for $(Y_n)$. There exist positive constants $c_{12}$, $c_{13}$, $c_{14}$, $c_{15}$, $c_{16}$ and $c_{17}$ such that for all $n\ge 1$,

 {\rm (i)} $\E(m^{Y_n}N_n^{(0)})\le c_{12} n$;
 
 {\rm (ii)} $c_{13} n^2\le \E(Y_n m^{Y_n}N_n^{(0)})\le c_{14} n^2$;

 {\rm (iii)} $\E(Y_n^2m^{Y_n}N_n^{(0)})\le c_{15} n^3$;
 
 {\rm (iv)} $\E( m^{Y_n}(N_n^{(0)})^2)\le c_{16} n^3$;
 
 {\rm (v)} $\E( Y_n m^{Y_n}(N_n^{(0)})^2)\le c_{17} n^4$.

\end{lemma}

\medskip

\noindent {\it Proof.} (i) Let $\mu_n(k) := m^{-k}\, \E([1-(m-1)Y_n]m^{Y_n} \, {\bf 1}_{\{X_n\le k\}})$ as in \eqref{mu}. Since the system is critical, we have $(m-1)\E(Y_n\, m^{Y_n}) =\E(m^{Y_n})$, so $\mu_n(k) = m^{-k}\E [((m-1)Y_n-1)m^{Y_n}\, {\bf 1}_{\{ Y_n\ge k+1\} }]$.
  
We have $\E(m^{Y_n}N_n^{(0)}) = \sum_{k=0}^\infty m^k \, \E(N_n^{(0)}\, {\bf 1}_{\{Y_n=k\}})$ by the Fubini--Tonelli theorem. Recall from \eqref{rel_pivot_openpaths} that $\E(N_n^{(0)} \, {\bf 1}_{\{ Y_n =k\} }) = \mu_0(0) \, \E(S_n^{(0,k)})$, which, by Theorem \ref{t:new_iteration}, is bounded by $\mu_n(k)$. Hence
\begin{eqnarray*}
    \E(m^{Y_n}N_n^{(0)}) 
 &\le& \sum_{k=0}^\infty \E [((m-1)Y_n-1)m^{Y_n}\, {\bf 1}_{\{ Y_n\ge k+1\} }]
    \\
 &=& \E [((m-1)Y_n-1)m^{Y_n}\, \sum_{k=0}^\infty {\bf 1}_{\{ Y_n\ge k+1\} }]
    \\
 &=& \E [Y_n((m-1)Y_n-1)m^{Y_n} \, {\bf 1}_{\{ Y_n \ge 1\} }]
    \\
 &\le& (m-1) \E (Y_n^2 m^{Y_n})\, .
\end{eqnarray*}

\noindent By Lemma \ref{l:G_nk123}~(iii), we have $\E(Y_n^2 m^{Y_n})\le c_4\, n$. Thus $\E(m^{Y_n}N_n^{(0)}) \le (m-1)c_4\, n$.

(ii)--(v) The proofs of (ii), (iii), (iv) and (v) are similar, all originating from \eqref{pf:recursive_formula3} (which we apply to the recursive system $(Y_n)$ in lieu of $(X_n)$). Let $i\ge 0$ be an integer, and write $G_n(s, \, t) := \E(s^{Y_n} \, t^{N^{(i)}_n})$. Define
$$
I_n(s)
    :=
    \E(s^{Y_n} \, N_n^{(i)})
    =
    \frac{\partial G_n(s, \, t)}{\partial t} \Big| \, _{t=1}  \, .
$$

\noindent It follows from \eqref{pf:recursive_formula3} that
\begin{equation}
    I_{n+1}(s)
    =
    \frac{m}{s} \, I_n(s) G_n(s)^{m-1}
    -
    \frac{m}{s} \, I_n(0) G_n(0)^{m-1} ,
    \label{e:In_recursion}
\end{equation}

\noindent where $G_n(s) := \E(s^{Y_n}) = G_n(s, \, 1)$. Multiplying by $s$ on both sides of \eqref{e:In_recursion}, and differentiating with respect to $s$, we get
$$
    I_{n+1}(s) + s \, I_{n+1}'(s)
    =
    m(m-1) G_n(s)^{m-2} G_n'(s) I_n(s)
    +
    m G_n(s)^{m-1} I_n'(s).
$$

\noindent At criticality, we have $m(m-1)G_n'(m)=G_n(m)$, which yields, for $n\ge 0$,
$$
    I_{n+1}(m) + mI_{n+1}'(m)
    =
    G_n(m)^{m-1}(I_n(m) + mI_n'(m)) \, .
$$

\noindent Iterating the identity, and noting that $\E[m^{Y_n}(1+Y_n)N_n^{(i)}] = I_n(m)+mI_n'(m)$, we get, for all $i\ge 0$ and $n\ge 1$,
\begin{eqnarray}
    \E[m^{Y_n}(1+Y_n)N_n^{(i)}]
 &=& \E[m^{Y_0}(1+Y_0)N_0^{(i)}] \prod_{j=0}^{n-1} G_j(m)^{m-1}
    \nonumber
    \\
 &=& (i+1)m^i \, \P(Y_0=i) \prod_{j=0}^{n-1} G_j(m)^{m-1} \, .
    \label{recursion_E(YN)}
\end{eqnarray}

\noindent This formula is known (Proposition 5.1 in \cite{bmvxyz_conjecture_DR}). We take $i=0$ from now on. 

Since $c_7 \, n^2\le \prod_{j=0}^{n-1}G_j(m)^{m-1}\le c_8 \, n^2$ (by Lemma \ref{l:G_nk123}~(v)), (ii) follows immediately. 

To prove (iii), we note from \eqref{e:In_recursion} that
$$
\frac{\partial^2 (sI_{n+1}(s))}{\partial s^2}
=
m \, \frac{\partial^2}{\partial s^2} [G_n(s)^{m-1}I_n(s)].
$$

\noindent Taking $s=m$, we get
\begin{eqnarray*}
    2I_{n+1}'(m)+mI_{n+1}''(m)
 &=& mG_n(m)^{m-1}I_n''(m)+2m(m-1)G_n(m)^{m-2}G_n'(m)I_n'(m) 
    \\
 && +
    m(m-1)G_n(m)^{m-2}G_n''(m)I_n(m)
    \\
 && +
    m(m-1)(m-2)G_n(m)^{m-3}G_n'(m)^2I_n(m).
\end{eqnarray*}

\noindent On the left-hand side, the expression is at least $mI_{n+1}''(m)$, whereas on the right-hand side, for the third term, we use $G_n(m)^{m-2} \le G_n(m)^{m-1}$. Since $G_n(m)=m(m-1)G_n'(m)$ at criticality, this yields that
\begin{eqnarray*}
    I_{n+1}''(m)
 &\le& G_n(m)^{m-1}I_n''(m)
    \\
 && +
    G_n(m)^{m-1} \Big( \frac{2}{m} \, I_n'(m)+ (m-1)I_n(m)G_n''(m)+ \frac{m-2}{m^2(m-1)} \, I_n(m) \Big).
\end{eqnarray*}

\noindent For the second expression on the right-hand side, we have, by Lemma \ref{l:G_nk123}~(i), $G_n(m)^{m-1} \le c_2^{m-1}$, and by (ii) which has just been proved, $I_n'(m)=\E(Y_n m^{Y_n-1} N_n^{(0)})\le \frac{c_{14}}{m} \, n^2$, whereas by (i), $I_n(m)=\E(m^{Y_n}N_n^{(0)})\le c_{12} \, n$, and by Lemma \ref{l:G_nk123}~(iii), $G_n''(m)\le \E(Y_n^2m^{Y_n})\le c_4 \, n$. Consequently,  
\begin{equation}
    I_{n+1}''(m)
    \le
    G_n(m)^{m-1} I_n''(m) + c_{18} \, n^2,
    \label{pf:Lemma3.2(iii)}
\end{equation}

\noindent where $c_{18} := c_2^{m-1} (\frac{2}{m^2}c_{14} + (m-1) c_{12} c_4 + \frac{m-2}{m^2(m-1)}c_{12})$. Iterating the inequality yields that for $n\ge 1$ (with $\prod_\varnothing := 1$),
$$
I_n''(m)
\le 
I_1''(m) \prod_{k=1}^{n-1} G_k(m)^{m-1}
+ 
c_{18}\sum_{\ell=1}^{n-1} \ell^2 \prod_{k=\ell+1}^{n-1}G_k(m)^{m-1} .
$$

\noindent By Lemma \ref{l:G_nk123}~(v), $\prod_{k=\ell+1}^{n-1}G_k(m)^{m-1}\le \frac{c_8n^2}{c_7(\ell+1)^2}$, this gives that for $n\ge 1$,
$$
I_n''(m)
\le
I_1''(m) c_8 n^2+ \frac{c_{18}c_8}{c_7} \, \sum_{\ell=1}^{n-1} \ell^2 \frac{n^2}{(\ell+1)^2}
\le
c_8 \big(\frac{c_{18}}{c_7}+I_1''(m)\big) \, n^3.
$$

\noindent Together with (ii) and (i), this implies (iii).

To prove (iv), let us define
$$
J_n(s)
:= 
\E(s^{Y_n} (N_n^{(0)})^2)
=
\frac{\partial^2}{\partial t^2} \, G_n(s, \, t) \, \Big| \, _{t=1}
+
I_n(s) \, .
$$

\noindent By \eqref{pf:recursive_formula3},
\begin{eqnarray*}
    J_{n+1}(m)
 &=& G_n(m)^{m-1} J_n(m)
    +
    (m-1) \, G_n(m)^{m-2} I_n(m)^2
    \\
 &&\qquad
    -
    (m-1) \, G_n(0)^{m-2} I_n(0)^2
    -
    G_n(0)^{m-1} J_n(0) 
    \\
 &\le& G_n(m)^{m-1} J_n(m)
    +
    (m-1) \, G_n(m)^{m-2} I_n(m)^2 \, .
\end{eqnarray*}

\noindent For the last term on the right-hand side, we use $G_n(m)^{m-2} \le c_2^{m-2}$ (by Lemma \ref{l:G_nk123}~(i)), and $I_n(m)=\E(m^{Y_n}N_n^{(0)})\le c_{12} \, n$ (by (i)). Therefore, with $c_{19} := (m-1) c_2^{m-2} c_{12}^2$,
$$
J_{n+1}(m)
\le
G_n(m)^{m-1} J_n(m) + c_{19} \, n^2 \, ,
$$

\noindent which is the analogue of \eqref{pf:Lemma3.2(iii)} for $J_n(m)$ in place of $I_n''(m)$. Since $J_0(m)=\E(m^{Y_0}(N_0^{(0)})^2) = \P(Y_0=0)$, this yields $\E(m^{Y_n}(N_n^{(0)})^2) = J_n(m) \le c_{16}\, n^3$ for some $c_{16}>0$ and all $n\ge 1$, exactly as \eqref{pf:Lemma3.2(iii)} has implied (iii).

It remains to check (v), which is similar, with more computations. Define
$$
K_n(s) 
:=
\E(Y_n s^{Y_n} (N_n^{(0)})^2)
=
s\frac{\partial}{\partial s} \frac{\partial^2}{\partial t^2} \, G_n(s, \, t) \, \Big| \, _{t=1}
+
s I_n'(s) \, .
$$

\noindent By \eqref{pf:recursive_formula3} and the identity $G_n'(m) = \frac{G_n(m)}{m(m-1)}$ at criticality, we get  \begin{eqnarray*}
    K_{n+1}(m)
 &=& G_n(m)^{m-1} K_n(m) 
    - 
    I_n(m)^2G_n(m)^{m-2} 
    + 
    2m(m-1) I_n(m)I_n'(m)G_n(m)^{m-2}
    \\
 && + 
    J_n(0)G_n(0)^{m-1} 
    + 
    (m-1)I_n(0)^2G_n(0)^{m-2} \, .
\end{eqnarray*}

\noindent On the right-hand side, we throw away the second term due to its negativity, use for the third term $I_n(m) \le c_{12} \, n$ (by (i)), $I_n'(m) \le \frac{c_{14}}{m} \, n^2$ (by (ii)) and $G_n(m)^{m-2} \le c_2^{m-2}$ (by Lemma \ref{l:G_nk123}~(i)), and for the last two terms, $G_n(0) = \P (Y_n=0) \le 1$, $J_n(0) \le J_n(m) \le c_{16} \, n^3$ (by (iv)), and $I_n(0)^2 \le I_n(m)^2 \le c_{12}^2 \, n^2$ (by (i) again). This leads to, with $c_{20} := 2(m-1)c_{12} c_{14}c_2^{m-2} + c_{16} + (m-1)c_{12}^2$: for $n\ge 1$,
$$
K_{n+1}(m) \le G_n(m)^{m-1} K_n(m) + c_{20} \, n^3 \, ,
$$ 

\noindent which is once again the analogue of \eqref{pf:Lemma3.2(iii)} for $K_n(m)$ in place of $I_n''(m)$. This yields that $\E(Y_n m^{Y_n}(N_n^{(0)})^2) = K_n(m) \le c_{21}\, n^4$ for some $c_{21}>0$ and all $n\ge 1$ exactly as \eqref{pf:Lemma3.2(iii)} has implied (iii). The proof of Lemma \ref{l:YLNUPP} is complete.\qed

\bigskip

We have yet two more lemmas, Lemmas \ref{l:moment_derivative_lb} and \ref{l:truncated_moment_lb} to present in this section. Our interest is in Lemma \ref{l:truncated_moment_lb}, which (in the case $k=1$) will be used a couple of times in Section \ref{s:lb} in the proof of the lower bound in Theorem \ref{t:main}. The proof of Lemma \ref{l:truncated_moment_lb} relies on Lemma \ref{l:moment_derivative_lb}. 

\medskip

\begin{lemma}
\label{l:moment_derivative_lb}

 Assume $\E(t^{Y_0})<\infty$ for some $t>m$. Let $G_n(s) := \E (s^{Y_n})$ as before.
 
 {\rm (i)} There exists a constant $c_{22}>0$ such that for all integers $n\ge 1$ and $k\ge 1$,
 $$
 G_n^{(k)}(m) \le k! \, n^{k-1} \, \ee^{c_{22} k} \, .
 $$
 
 {\rm (ii)} There exist constants $c_{23} \ge 1$ and $c_{24}>0$ such that for integers $k\ge 1$ and $n\ge \ee^{c_{23}k}$,
 $$
 G_n^{(k)}(m)
 \ge
 k! \, n^{k-1} \, \ee^{-c_{24} k} \, .
 $$

\end{lemma}

\medskip

\noindent {\it Proof.} The inequality in (i) was proved in \cite[Theorem 2]{xz_stable}. Let us prove (ii). For fixed $k$, (ii) is simply the first inequality in \eqref{Yn_moments}, so we need to prove it only for $k \ge 6$.

We claim the existence of constants $c_{23} \ge 1$ and $c_{25}>0$ such that for integers $k \ge 2$ and $n\ge \lceil \ee^{c_{23}k} \rceil$,  \begin{equation}
     G_n^{(k)}(m)
     \ge
     k! \, n^{k-1} \, \ee^{-c_{25} (k-1)} \, ,
     \label{moment_derivative_lb}
\end{equation}

\noindent while keeping in mind that \eqref{moment_derivative_lb} is already proved for $2\le k\le 5$. Our proof is based on an argument by induction in $k$. 

Let $k \ge 6$. Assume \eqref{moment_derivative_lb} is valid for all derivatives of order $\le k-1$ at $s=m$. We need to prove it for $G_n^{(k)}(m)$. Let us recall the formula \eqref{recurrence_derivatives_Gn} at $s=m$:
$$
m G_{n+1}^{(k)}(m) + k G_{n+1}^{(k-1)}(m)
=
\sum_{(k_1, \ldots, k_m) \in \z_+^m: \, k_1 + \cdots + k_m = k} \frac{k!}{k_1! \cdots k_m!} \prod_{i=1}^m G_n^{(k_i)}(m) \, .
$$

\noindent On the right-hand side, we only consider two special types of $(k_1, \ldots, k_m) \in \z_+^m$ satisfying $k_1 +\cdots+k_m=k$: (1) one of $k_i$ is $k$ (so all others vanish); (2) one of $k_i$ is $\ell$ (for some $1\le \ell \le \lfloor \frac{k}{2} \rfloor -1$), another is $k-\ell$ (and all others vanish). Hence, for $k \ge 6$,
\begin{eqnarray*}
    m G_{n+1}^{(k)}(m) + k G_{n+1}^{(k-1)}(m)
 &\ge& m G_n^{(k)}(m) G_n(m)^{m-1}
    \\
 && +
    m(m-1) \sum_{\ell=1}^{\lfloor \frac{k}{2}\rfloor -1} \frac{k!}{\ell! \, (k-\ell)!} \, G_n^{(\ell)}(m) G_n^{(k-\ell)}(m) G_n(m)^{m-2} \, .
\end{eqnarray*}

\noindent On the right-hand side, the term $\ell=1$ is special: we have $G_n'(m) = \frac{G_n(m)}{m(m-1)}$ at criticality. This yields
\begin{eqnarray*}
    m G_{n+1}^{(k)}(m) + k G_{n+1}^{(k-1)}(m)
 &\ge& [m G_n^{(k)}(m) + k G_n^{(k-1)}(m)] G_n(m)^{m-1}
    \\
 && +
    m(m-1) \sum_{\ell=2}^{\lfloor \frac{k}{2} \rfloor-1} \frac{k!}{\ell! \, (k-\ell)!} \, G_n^{(\ell)}(m) G_n^{(k-\ell)}(m) G_n(m)^{m-2} \, .
\end{eqnarray*}

\noindent Let $2\le \ell\le \lfloor \frac{k}{2} \rfloor -1$. By the induction assumption, if $n\ge \lceil \ee^{c_{23} (k-2)} \rceil$, then $G_n^{(\ell)}(m) G_n^{(k-\ell)}(m) \ge \ell!\, (k-\ell)! \, n^{k-2} \ee^{-c_{25} (k-2)}$. We also use the trivial inequalities $m(m-1) \ge 1$ and $G_n(m)^{m-2} \ge 1$. Consequently,
\begin{eqnarray*}
    m(m-1) \sum_{\ell=2}^{\lfloor \frac{k}{2} \rfloor -1} \frac{k!}{\ell! \, (k-\ell)!} \, G_n^{(\ell)}(m) G_n^{(k-\ell)}(m) G_n(m)^{m-2}
 &\ge& \sum_{\ell=2}^{\lfloor \frac{k}{2} \rfloor -1} k! \, n^{k-2} \ee^{-c_{25} (k-2)}
    \\
 &=& (\lfloor \frac{k}{2} \rfloor -2)\, k! \, n^{k-2} \ee^{-c_{25} (k-2)} \, .
\end{eqnarray*}

\noindent Since $k \ge 6$, we have $\lfloor \frac{k}{2} \rfloor -2  \ge \frac{k}{7}$. Hence, for $k \ge 6$ and $n\ge \lceil \ee^{c_{23} (k-2)} \rceil$,
$$
m G_{n+1}^{(k)}(m) + k G_{n+1}^{(k-1)}(m)
\ge
[m G_n^{(k)}(m) + k G_n^{(k-1)}(m)] G_n(m)^{m-1}
+
\frac{k! \, k}{7}\, n^{k-2} \ee^{-c_{25} (k-2)} \, .
$$

\noindent Iterating the inequality gives that, for $n\ge n_0+1$ where $n_0 = n_0(k) :=\lceil \ee^{c_{23} (k-2)} \rceil$,
\begin{eqnarray*}
 &&m G_n^{(k)}(m) + k G_n^{(k-1)}(m)
    \\
 &\ge& \Big( m G_{n_0}^{(k)}(m) + k G_{n_0}^{(k-1)}(m) + \sum_{i=n_0}^{n-1} \frac{\frac{k! \, k}{7} i^{k-2} \ee^{-c_{25} (k-2)}}{\prod_{j=n_0}^i G_j(m)^{m-1}} \Big) \prod_{j=n_0}^{n-1} G_j(m)^{m-1} 
    \\
 &\ge& \sum_{i=n_0}^{n-1} \frac{\frac{k! \, k}{7} i^{k-2} \ee^{-c_{25} (k-2)}}{\prod_{j=n_0}^i G_j(m)^{m-1}} \prod_{j=n_0}^{n-1} G_j(m)^{m-1} \, . 
\end{eqnarray*}

\noindent By Lemma \ref{l:G_nk123}~(v), $c_7 \, n^2\le \prod_{i=0}^{n-1}G_i(m)^{m-1}\le c_8 \, n^2$, so 
$$
m G_n^{(k)}(m) + k G_n^{(k-1)}(m)
\ge
\frac{k! \, k}{7} \ee^{-c_{25}(k-2)} \sum_{i=n_0}^{n-1} \frac{i^{k-2}}{(i+1)^2} \, \frac{c_7}{c_8} \, n^2 \, .
$$

\noindent If $n\ge \lceil \ee^{c_{23} k} \rceil$, then $n_0 \le \frac{n}{2}$ (recalling that $c_{23}\ge 1$), so $\sum_{i=n_0}^{n-1} \frac{i^{k-2}}{(i+1)^2} \ge \sum_{i= \lfloor n/2\rfloor}^{n-1} \frac{i^{k-2}}{(i+1)^2}  \ge c_{26} \, \frac{n^{k-3}}{k}$ for some constant $c_{26}>0$. This yields, for $n\ge \lceil \ee^{c_{23} k} \rceil$ (with $c_{27} := \frac{c_7 c_{26}}{7c_8}$)
\begin{equation}
    m G_n^{(k)}(m) + k G_n^{(k-1)}(m)
    \ge
    c_{27} \, k! \, \ee^{-c_{25}(k-2)} n^{k-1} \, .
    \label{pf:derivatives_Gn}
\end{equation}

\noindent On the other hand, by (i), $G_n^{(k-1)}(m) \le \ee^{c_{22} (k-1)} \, (k-1)! \, n^{k-2}$ for some constant $c_{22}>0$. 
Thus
$$
k G_n^{(k-1)}(m) \le k \, \ee^{c_{22} (k-1)} \, (k-1)! \, n^{k-2} ,
$$

\noindent which is bounded by $\frac{c_{27}}{2} \, k! \, \ee^{-c_{25} (k-2)} n^{k-1}$ if $n\ge \lceil \ee^{c_{23} k} \rceil$ (with a large constant $c_{23}$). Going back to \eqref{pf:derivatives_Gn}, we obtain: for $k \ge 6$ and $n\ge \lceil \ee^{c_{23} k} \rceil$,
$$
m G_n^{(k)}(m)
\ge
\frac{c_{27}}{2} \, k! \, \ee^{-c_{25} (k-2)} n^{k-1}\, ,
$$

\noindent which is at least $m \, k! \, \ee^{-c_{25} (k-1)} n^{k-1}$ if $c_{25}$ is chosen to be sufficiently large.\footnote{The constants are chosen in this order: first $c_{25}$ then $c_{23}$, while our presentation is in the opposite order.} As a consequence, \eqref{moment_derivative_lb} is valid for $k \ge 6$. 

Lemma \ref{l:moment_derivative_lb} is proved.\qed

\medskip

\begin{lemma}
\label{l:truncated_moment_lb}

 Assume $\E(t^{Y_0})<\infty$ for some $t>m$. For any $\lambda>0$ and any integer $k\ge 1$,
 \begin{equation}
      \liminf_{n\to \infty} \, \frac{1}{n^{k-1}} \, \E(Y_n^k m^{Y_n} \, {\bf 1}_{\{ Y_n \ge \lambda n\} })
      >
      0 \, .
      \label{truncated_moment_lb}
 \end{equation}

\end{lemma}

\medskip

\noindent {\it Proof.} We have, for $k\ge 1$,
$$
\E(Y_n^k m^{Y_n})
\le
\E(Y_n^k m^{Y_n} \, {\bf 1}_{\{ Y_n \ge \lambda n \} })
+
(\lambda n)^{k-1} \, \E(Y_n m^{Y_n})\, .
$$

\noindent By Lemma \ref{l:G_nk123}~(ii), $\E(Y_n \, m^{Y_n})\le c_3$. Thus
$$
\E(Y_n^k m^{Y_n})
\le
\E(Y_n^k m^{Y_n} \, {\bf 1}_{\{ Y_n \ge \lambda n \} })
+
c_3(\lambda n)^{k-1} \, .
$$

\noindent On the other hand, by Lemma \ref{l:moment_derivative_lb}~(ii), for all sufficiently large integer $k$ (say $k \ge k_0$) and all integer $n\ge \ee^{c_{23} k}$,
$$
\E(Y_n^k m^{Y_n}) \ge 2 c_3 (\lambda n)^{k-1} \, .
$$

\noindent Consequently, for $k \ge k_0$ and $n\ge \ee^{c_{23} k}$,
$$
2 c_3 (\lambda n)^{k-1}
\le
\E(Y_n^k m^{Y_n})
\le
\E(Y_n^k m^{Y_n} \, {\bf 1}_{\{ Y_n \ge \lambda n \} })
+
c_3 (\lambda n)^{k-1} \, ,
$$

\noindent which implies that
\begin{equation}
    \E(Y_n^k m^{Y_n} \, {\bf 1}_{\{ Y_n \ge \lambda n \} })
    \ge
    c_3 (\lambda n)^{k-1},
    \qquad
    k \ge k_0, \; n\ge \ee^{c_{23} k}\, .
    \label{pf_moment_derivative_lb}
\end{equation}

\noindent In particular,
\begin{equation}
    \liminf_{n\to \infty} \, \frac{1}{n^{k-1}} \, \E(Y_n^k m^{Y_n} \, {\bf 1}_{\{ Y_n \ge \lambda n\} })
    >
    0 \, ,
    \label{moment_derivative_lb_bis}
\end{equation}

\noindent which yields Lemma \ref{l:truncated_moment_lb} when $k \ge k_0$.

To remove the restriction $k \ge k_0$ for the validity of \eqref{moment_derivative_lb_bis}, we observe that by the Cauchy--Schwarz inequality, 
$$
\E(Y_n^k m^{Y_n} \, {\bf 1}_{\{ Y_n \ge \lambda n \} })
\le
[\E(Y_n^{k+1} m^{Y_n})]^{1/2} \,
[\E(Y_n^{k-1} m^{Y_n} \, {\bf 1}_{\{ Y_n \ge \lambda n \} })]^{1/2} \, .
$$

\noindent [A similar inequality was already used in \eqref{CS}.] Recall from \eqref{Yn_moments} that $\E(Y_n^{k+1} m^{Y_n}) \le c_{10}\, n^k$. Hence, for $k \ge 1$ and $n\ge 1$,
\begin{equation}
    \E(Y_n^{k-1} m^{Y_n} \, {\bf 1}_{\{ Y_n \ge \lambda n \} })
    \ge
    \frac{[\E(Y_n^k m^{Y_n} \, {\bf 1}_{\{ Y_n \ge \lambda n \} })]^2}{c_{10}\, n^k} \, .
    \label{pf_moment_derivative_lb2}
\end{equation}

\noindent Consequently, the inequality \eqref{moment_derivative_lb_bis} applied to $k =k_0$, together with \eqref{pf_moment_derivative_lb2}, yield that \eqref{moment_derivative_lb_bis} is valid for $k = k_0-1$. Iterating the argument, we see that \eqref{moment_derivative_lb_bis} remains valid for all $k \ge 1$.\qed

\section{Proof of Theorem \ref{t:main}: upper bound}
\label{s:ub}

Let $(Y_n, \, n\ge 0)$ be a Derrida--Retaux system satisfying $\E(m^{Y_0}) = (m-1) \E(Y_0\, m^{Y_0})<\infty$; so it is critical. The upper bound in Theorem \ref{t:main} is a consequence of the following result in the particular case $\ell=0$.

\medskip

\begin{proposition}
 \label{p:ub}

 Assume $\E(Y_0^5 \, m^{Y_0})<\infty$. There exists $c_{28}>0$ such that for integers $n\ge 2$ and $\ell\ge 0$, 
 \begin{equation}
     \P(Y_n \ge \ell+1)
     \le
     \frac{c_{28}}{m^\ell} \, \frac{(\log n)^2}{n^2} \, .
     \label{p:ub1}
 \end{equation}

\end{proposition}

\medskip

In the proof of Proposition \ref{p:ub}, we are going to truncate the initial distribution of a system and compare it with the system itself. Let $(X_n, \, n\ge 0)$ be a Derrida--Retaux system, with genealogical tree $\T$. Let $B \subset \{ 0, \, 1, \, 2, \, \ldots\}$. Define
$$
\widehat{X}_0 := X_0 \, {\bf 1}_{\{ X_0 \notin B\} } \, .
$$

\noindent Let $(\widehat{X}_n, \, n\ge 0)$ be a Derrida--Retaux system with initial distribution $\widehat{X}_0$, associated with genealogical tree $\T$, such that $X(v) \ge \widehat{X}(v)$ for all $v\in \T$. In particular, $\widehat{X}(v) = X(v) \, {\bf 1}_{\{ X(v) \notin B\} }$ if $|v|=0$.

\medskip

\begin{lemma}
\label{l:coupling}

 Fix $n\ge 1$. If there are no open paths in the system $(X_k, \, k\ge 0)$ that lead to vertex $\mathfrak{e}_n$ and start from a vertex at the initial generation with a value lying in $B$, then $X(\mathfrak{e}_n) = \widehat{X}(\mathfrak{e}_n)$.

\end{lemma}

\medskip

\noindent {\it Proof of Lemma \ref{l:coupling}.} Consider the system $(X_k, \, k\ge 0)$. Let $R_n$ denote the set of vertices (including $\mathfrak{e}_n$ itself) in all the open paths starting from the initial generation and leading to $\mathfrak{e}_n$. Let $R_{n,0}$ be the set of elements of $R_n$ that are in the initial generation. [So $N_n$, defined in \eqref{N}, is simply $|R_{n,0}|$, the cardinality of $R_{n,0}$.] Then the value of $X(\mathfrak{e}_n)$ does not depend on $\{ X(v), \, v\in \T_n\backslash R_{n,0} \hbox{ \rm with }|v|=0\}$, and
\begin{equation}
    X(\mathfrak{e}_n) = \sum_{v\in R_{n,0}} X(v)  - |R_n \backslash R_{n,0}| \, .
    \label{X_representation}
\end{equation}

\noindent [If $R_{n,0}$ is empty, then so is $R_n$, and \eqref{X_representation} holds trivially.]


Assume that there are no open paths in the system $(X_k, \, k\ge 0)$ that lead to vertex $\mathfrak{e}_n$ and start from a vertex at the initial generation with a value lying in $B$. Then for any $v\in R_{n,0}$, we have $X(v) \notin B$, so that $\widehat{X}(v) = X(v)$. This implies that $\widehat{X}(v) \ge X(v)$ for $v\in R_n$. In particular, $\widehat{X}(\mathfrak{e}_n) \ge X(\mathfrak{e}_n)$.\qed

\bigskip

In order to prove Proposition \ref{p:ub}, we work with a truncated version of the system $(Y_n, \, n\ge 0)$, the level of the truncation depending on an integer parameter $M$ which is sufficiently large (we write generically $M\ge M_0$; the value of $M_0$ may change from line to line). For integer $M\ge M_0$, define
\begin{equation}
    \zeta(M)
    :=
    -\log \E(Y_0^3m^{Y_0} \, {\bf 1}_{\{ Y_0>M \}}) 
    \in [1, \, \infty]\, .
    \label{zeta}
\end{equation}

\noindent [Notation: $\zeta(M) := \infty$ if $\P(Y_0>M)=0$. We have $\lim_{M\to \infty} \zeta(M) = \infty$ under the assumption $\E(Y_0^3m^{Y_0})<\infty$.] Define
\begin{equation}
    Z_0 = Z_0(M) := Y_0\, {\bf 1}_{\{ Y_0\le M\zeta(M) \}} \, .
    \label{Z0}
\end{equation}

\noindent We consider a new Derrida--Retaux system $(Z_n, \, n\ge 0)$ with initial distribution $Z_0$. Since $Z_0 \le Y_0$ and $(Y_n, \, n\ge 0)$ is critical, the new system $(Z_n, \, n\ge 0)$ is subcritical, or critical if $Y_0\le M\zeta(M)$ a.s.

\bigskip

\noindent {\it Outline of the proof of Proposition \ref{p:ub}.} The key ingredient in the proof is the following estimate (Lemma \ref{l:P(Z>k)}): if $\E(Y_0^3m^{Y_0})<\infty$, there exist constants $c_{43}>0$, $c_{44}>0$ and $M_0\ge 1$ such that for all integers $n\ge M\ge M_0$,
$$
    \P(Z_n \ge \ell+1)
    \le 
    c_{43} \Big( \frac{1}{M^2} + \frac{1}{M}\, \ee^{- c_{44}\frac{n}{M}} \Big) m^{-\ell},
    \qquad
    \ell\ge 0\, .
    \leqno{\hbox{\it(\ref{P(Z>k)})}}
$$

\noindent This will yield \eqref{p:ub1} by taking $M$ to be (the integer part of) a suitable constant multiple of $\frac{n}{\log n}$, because we will see that the difference between $\P(Y_n \ge \ell+1)$ and $\P(Z_n \ge \ell+1)$ can be controlled under further integrability assumption $\E(Y_0^5m^{Y_0})<\infty$. 

To prove \eqref{P(Z>k)}, we define a deterministic set $B_M$ of integers $y$ of order of magnitude $M$ such that conditionally on $Z_M=y$, ``with reasonable probability", there are at least (a constant multiple of) $M^2$ open paths from the initial generation to generation $M$. 

We distinguish two possible situations:

\medskip

{\leftskip=1.5truecm \rightskip=1.5truecm

\noindent Case 1: there are no open paths starting from generation $M$ 

\noindent \phantom{Case 1: }with a value lying in $B_M$ leading to generation $n+M$;

\noindent Case 2: there are such open paths.

}

\medskip

It is not hard to treat Case 2. The basic idea is that if there are such open paths, then by the construction of $B_M$, with ``reasonable probability", there are at least (a constant multiple of) $M^2$ open paths leading to generation $n+M$ and starting from the initial generation with value $0$. We are going to see that thanks to \eqref{E(N_Y>k)}, this leads to:
$$
\P(Z_{n+M} \ge \ell, \; \hbox{Case 2})
\le
\frac{c}{M^2} \, \frac{1}{m^\ell} \, ,
$$

\noindent for some constant $c>0$, all $n\ge 0$ and $\ell \ge 1$. 

To treat Case 1 which is more delicate, we define
$$
W_0 = W_0(M) := Z_M \, {\bf 1}_{\{ Z_M \notin B_M\} } \, ,
$$

\noindent and consider yet another Derrida--Retaux system $(W_n, \, n\ge 0)$ with initial distribution $W_0$. The reason for which we are interested in this new system is that by Lemma \ref{l:coupling}, in Case 1, if $Z_{n+M} =\ell$, then $W_n =\ell$, so that for all $n\ge 0$,
$$
\P(Z_{n+M} =\ell, \; \hbox{Case 1})
\le
\P(W_n =\ell)
\le
\frac{\E(m^{W_n})-1}{m^\ell -1} \, .
$$

\noindent We need to control $\E(m^{W_n})$. This relies on subcriticality of the $(W_i, \, i\ge 0)$, and on a simple estimate (Lemma \ref{l:subcritcalspeed}) on the moment generating function for all subcritical Derrida--Retaux systems satisfying a convenient condition for the initial distribution. We check by means of a technical lemma (Lemma \ref{l:gen_fct_W}) that $W_{n_0}$ satisfies this condition, where $n_0 := \lfloor c_{34} M \rfloor$ for some constant $c_{34}>0$ and all sufficiently large $M$.\footnote{It is not a surprise to consider the system at the intermediate generation $n_0$. A similar idea is going to be used in the proof of the lower bound in Section \ref{s:lb}.} This will yield that there exist constants $c_{33}>0$ and $c_{35}\in (0, \, 1)$, such that for all integer $n\ge 0$,
$$
\E(m^{W_{n+n_0}}) \le 1+ \frac{c_{33}}{M} \Big( 1- \frac{c_{35}}{M} \Big)^n \, .
$$

\noindent Consequently, for $\ell \ge 1$ and $n\ge 2n_0$,
$$
\P(Z_{n+M} =\ell, \; \hbox{Case 1})
\le
\frac{\E(m^{W_n})-1}{m^\ell -1}
\le
\frac{\frac{c_{33}}{M} ( 1- \frac{c_{35}}{M})^{n-n_0}}{m^\ell -1}
\le
\frac{\frac{c_{33}}{M} ( 1- \frac{c_{35}}{M})^{n/2}}{m^\ell -1} \, .
$$

\noindent Combining this with the bound for $\P(Z_{n+M} \ge \ell, \; \hbox{Case 2})
$, we see that for all $n\ge 2n_0$ and $\ell \ge 1$, 
$$
\P(Z_{n+M} \ge \ell)
\le
\frac{c}{M^2} \, \frac{1}{m^\ell}
+
\sum_{j=\ell}^\infty \frac{\frac{c_{33}}{M} ( 1- \frac{c_{35}}{M})^{n/2}}{m^j -1} \, .
$$

\noindent This will lead to \eqref{P(Z>k)}.\qed 


\bigskip

The rest of the section is devoted to the proof of Proposition \ref{p:ub}. For the sake of clarity, it is divided into three steps and presented in distinct parts. The first step is about localization; it consists in understanding where the truncated Derrida--Retaux system $(Z_n, \, n\ge 0)$ lies with good (weighted) probability conditioning on sustainability; this is then used to prove Lemma \ref{l:gen_fct_W}, the main result in the step. The second step gives \eqref{P(Z>k)} (which is Lemma~\ref{l:P(Z>k)}), whereas the third step completes the proof of Proposition \ref{p:ub}.

\subsection{Step 1. Localization}
\label{subs:localization}

The aim of Step 1 is to prove Lemma \ref{l:gen_fct_W}, which is the technical part of preparation to give an upper bound for $\P(Z_n \ge \ell+1)$ in Lemma~\ref{l:P(Z>k)}. 
For $n\ge 1$, let $N_n^{(0)}$ denote the number of open paths for the critical Derrida--Retaux system $(Y_n)$ (replacing $(X_n)$ by $(Y_n)$ in \eqref{N}), and let $L_n^{(0)}$ be the corresponding quantity for the truncated Derrida--Retaux system $(Z_n)$; define the deterministic set $B_n$ by
\begin{equation}
    B_n
    :=
    \{ y \in [c_{29}n, \, c_{30}n] \cap \z: \, \P(Z_n=y) >0, \; \P( L_n^{(0)}\ge c_{31} n^2 \, | \, Z_n =y)\ge c_{32}\} \, ,
    \label{BM}
\end{equation}

\noindent where the positive constants $c_{29}$, $c_{30}$, $c_{31}$ and $c_{32}$ are explicit but uninteresting: $c_{29} := \min\{ \frac{c_{13}}{5c_{12}}, \, \frac{5c_{15}}{2c_{13}}\}$, $c_{30} := \frac{5c_{15}}{c_{13}}$, $c_{31} := \min\{ \frac{c_{13}}{5c_3}, \, \frac{5c_{17}}{2c_{13}}\}$, $c_{32} := \frac{c_{13}^3}{250 c_3 c_{15} c_{17}} \min\{ \frac{c_{13}}{5c_{12}}, \, \frac{5c_{15}}{2c_{13}}\}$ ($c_3$ being from Lemma \ref{l:G_nk123}~(ii), and $c_{12}$, $c_{13}$, $c_{15}$ and $c_{17}$ from Lemma \ref{l:YLNUPP}). Let
\begin{equation}
    W_0 = W_0(M) := Z_M  \, {\bf 1}_{\{ Z_M \notin B_M\} } \, .
    \label{W0}
\end{equation}

\noindent We also consider a new Derrida--Retaux system $(W_n, \, n\ge 0)$ with initial distribution $W_0$. Since $W_0 \le Z_M$, and $(Z_n)$ is subcritical or critical, it follows that $(W_n, \, n\ge 0)$ is subcritical or critical. As such, we have three Derrida--Retaux systems: the critical system $(Y_n, \, n\ge 0)$ which we are interested in, and the two associated systems $(Z_n, \, n\ge 0)$ and $(W_n, \, n\ge 0)$. They are related in the following way: for all $n\ge 0$,
$$
W_n \le Z_{n+M} \le Y_{n+M} \, .
$$

\noindent We mention that both systems $(Z_n, \, n\ge 0)$ and $(W_n, \, n\ge 0)$ depend on $M$.

The aim of the step is to prove the following estimate. 

\medskip

\begin{lemma}
 \label{l:gen_fct_W} 
 
 Assume $\E(Y_0^3m^{Y_0})<\infty$. There exist constants $c_{33}\in (0, \, 1)$, $c_{34}>0$, $c_{35}\in (0, \, 1)$, such that for all sufficiently large integer $M$, with $s = s(M) := m+ \frac{c_{33}}{M}$ and $n_0 =n_0(M) := \lfloor c_{34} M \rfloor$,
 \begin{equation}
     \frac{m}{s} \, [\E(s^{W_{n_0}})]^{m-1}
     \le 
     1- \frac{c_{35}}{M} \, .
     \label{eq_l:gen_fct_W}
 \end{equation}

\end{lemma}

\medskip

The proof of Lemma \ref{l:gen_fct_W} relies on the following estimate.
\medskip

\begin{lemma}
 \label{l:B_nN_n>M^2} 
 
 Assume $\E(Y_0^3\, m^{Y_0})<\infty$.
 There exists a constant $c_{36}>0$ such that for all sufficiently large integer $M$,
 $$
 \E(m^{Z_M} \, {\bf 1}_{\{ Z_M\in B_M\} })
 \ge 
 \frac{c_{36}}{M}, 
 $$

 \noindent where $B_M$ is the deterministic set of integers defined in \eqref{BM}.\footnote{Since $m^{Z_M}\, {\bf 1}_{\{ Z_M \in B_M\} } \le m^{Z_M}\, {\bf 1}_{\{ Z_M \ge 1\} } \le m^{Y_M}\, {\bf 1}_{\{ Y_M \ge 1\} }$, it follows from Lemmas \ref{l:B_nN_n>M^2} and \ref{l:G_nk123}~(iv) that $\E(m^{Z_M}\, {\bf 1}_{\{ Z_M \ge 1\} })$ is of order of magnitude $\frac1M$. So the meaning of Lemma \ref{l:B_nN_n>M^2} is as follows: if we work with the probability measure having a density proportional to $m^{Z_M}$, then conditionally on the event $Z_M\ge 1$, with probability greater than a certain constant (independent of $M$), $Z_M$ lies in $B_M$.}
\end{lemma}

\medskip

\noindent {\it Proof of Lemma \ref{l:B_nN_n>M^2}.} Let $b_1$, $b_2$, $b_3$, $b_4$, $b_5$ be positive constants whose values will be chosen later. 
Let us consider the quantity
$$
q_n
:=
\E( Z_n m^{Z_n} L_n^{(0)} \, {\bf 1}_{\{ L_n^{(0)} \in [b_1 n^2, \, b_2 n^2], \, Z_n\in [b_3 n, \, b_4 n]\} } )\, .
$$

\noindent Assume $M\ge M_0$ ($M_0$ being a sufficiently large integer such that $\zeta(M)\ge 2$ for $M\ge M_0$) and $1\le n\le M$. We first recall that $Z_0=Y_0  \, {\bf 1}_{\{Y_0\le M\zeta(M)\}}$. Since $Y_n \ge Y_0-n \ge Y_0-M$ (this is why we assume $n\le M$), we have, on the event $\{ Y_n\le \frac12 M\zeta(M)\}$, $Y_0 \le Y_n + M \le M\zeta(M)$, so $Z_n=Y_n$ 
and $L_n^{(0)} = N_n^{(0)}$. Consequently, if $b_4 \le \frac{\zeta(M)}{2}$, then $\{ Y_n\in [b_3 n, \, b_4 n] \} \subset \{ Z_n\in [b_3 n, \, b_4 n]\}$, and 
\begin{eqnarray*}
    q_n 
 &\ge& \E( Y_n m^{Y_n} N_n^{(0)} \, {\bf 1}_{\{ N_n^{(0)} \in [b_1 n^2, \, b_2 n^2], \, Y_n\in [b_3 n, \, b_4 n]\} })
    \\
 &\ge& \E( Y_n m^{Y_n} N_n^{(0)}) 
    - 
    \E( Y_n m^{Y_n} N_n^{(0)}\, {\bf 1}_{\{ N_n^{(0)} \notin [b_1 n^2, \, b_2 n^2] \} } )
    -
    \E( Y_n m^{Y_n} N_n^{(0)}\, {\bf 1}_{\{ Y_n \notin [b_3 n, \, b_4 n] \} } )
    \\
 &=:& q_{n,1} - q_{n,2} - q_{n,3} \, ,
\end{eqnarray*}

\noindent with obvious notation. On the right-hand side, let us look at each of the three expressions $q_{n,i}$, $i\in \{1,\, 2, \, 3\}$. The first term $q_{n,1}$ is dealt with by Lemma \ref{l:YLNUPP}~(ii), which says that 
$$
q_{n,1} 
\ge
c_{13} n^2\, .
$$

\noindent For the second term $q_{n,2}$, we note that on the event $\{ N_n^{(0)} \notin [b_1 n^2, \, b_2 n^2] \}$, we have either $N_n^{(0)} < b_1 n^2$, in which case $Y_n m^{Y_n} N_n^{(0)} \le b_1 n^2 \, Y_n m^{Y_n}$, or $N_n^{(0)} > b_2 n^2$, in which case $Y_n m^{Y_n} N_n^{(0)} \le \frac{1}{b_2 n^2} \, Y_n m^{Y_n} (N_n^{(0)})^2$. Accordingly,
$$
q_{n,2}
\le
b_1 n^2 \, \E(Y_n m^{Y_n}) 
+
\frac{1}{b_2 n^2} \, \E(Y_n m^{Y_n} (N_n^{(0)})^2)
\le
c_3 b_1 n^2
+
\frac{c_{17}n^4}{b_2 n^2} ,
$$

\noindent using Lemma \ref{l:G_nk123}~(ii) and Lemma \ref{l:YLNUPP}~(v) for the last inequality. The third term $q_{n,3}$ is handled in the same way:
$$
q_{n,3}
\le
b_3 n \, \E(m^{Y_n} N_n^{(0)}) 
+
\frac{1}{b_4 n} \, \E(Y_n^2 m^{Y_n} N_n^{(0)})
\le
b_3 n \, c_{12} n
+
\frac{c_{15} n^3}{b_4 n} ,
$$

\noindent the last inequality following from Lemma \ref{l:YLNUPP}~(i) and (iii). Assembling these pieces yields that for $1\le n\le M$,
$$
q_n
\ge 
\Big( c_{13} - c_3 b_1 - \frac{c_{17}}{b_2} - c_{12} b_3 - \frac{c_{15}}{b_4} \Big) n^2
\ge
\frac{c_{13}}{5} \, n^2 \, ,
$$

\noindent if we choose $b_2 := \frac{5c_{17}}{c_{13}}$, $b_1 := \min\{ \frac{c_{13}}{5c_3}, \, \frac{b_2}{2}\}$, $b_4 := \frac{5c_{15}}{c_{13}}$ and $b_3 := \min\{ \frac{c_{13}}{5c_{12}}, \, \frac{b_4}{2} \}$ [$M$ being large enough so that $b_4 \le \frac{\zeta(M)}{2}$].

Let us have another look at $q_n := \E( Z_n m^{Z_n} L_n^{(0)} \, {\bf 1}_{\{ L_n^{(0)} \in [b_1 n^2, \, b_2 n^2], \, Z_n\in [b_3 n, \, b_4 n]\} } )$. On the event $\{ L_n^{(0)} \in [b_1 n^2, \, b_2 n^2], \, Z_n\in [b_3 n, \, b_4 n]\}$, we have $Z_n m^{Z_n} L_n^{(0)}\le b_2b_4n^3 \, m^{Z_n}$, so
$$
q_n 
\le 
b_2b_4n^3 \, \E( m^{Z_n} \, {\bf 1}_{\{ L_n^{(0)} \in [b_1 n^2, \, b_2 n^2], \, Z_n\in [b_3 n, \, b_4 n]\} } ) \, .
$$

\noindent Since $q_n \ge \frac{c_{13}}{5} \, n^2$, this yields that for $1\le n\le M$,
$$
\E( m^{Z_n} \, {\bf 1}_{\{ L_n^{(0)} \in [b_1 n^2, \, b_2 n^2], \, Z_n\in [b_3 n, \, b_4 n]\} } )
\ge 
\frac{c_{37}}{n} \, ,
$$

\noindent where $c_{37} := \frac{c_{13}}{5b_2b_4}$. A fortiori,
$$
\E(m^{Z_n} \, {\bf 1}_{\{ L_n^{(0)} \ge b_1 n^2, \, Z_n\in [b_3 n, \, b_4 n]\} })
\ge 
\frac{c_{37}}{n} \, .
$$

On the other hand,
\begin{eqnarray*}
 &&\E(m^{Z_n} \, {\bf 1}_{\{ L_n^{(0)} \ge b_1 n^2, \, Z_n\in [b_3 n, \, b_4 n]\} })
    \\
 &=& \E[ m^{Z_n} \, {\bf 1}_{\{ Z_n\in [b_3 n, \, b_4 n]\}}\, \P(L_n^{(0)} \ge b_1 n^2 \, | \, Z_n)]
    \\
 &\le& \E( m^{Z_n} \, {\bf 1}_{\{ Z_n\in [b_3 n, \, b_4 n]\}} \, {\bf 1}_{\{ \P( L_n^{(0)} \ge b_1 n^2 \, | \, Z_n)\ge b_5\}})
    +
    b_5 \, \E( m^{Z_n} \, {\bf 1}_{\{ Z_n\in [b_3 n, \, b_4 n]\}})
    \\
 &\le& \E( m^{Z_n} \, {\bf 1}_{\{ Z_n\in [b_3 n, \, b_4 n]\}} \, {\bf 1}_{\{ \P( L_n^{(0)} \ge b_1 n^2 \, | \, Z_n) \ge b_5\}})
    +
    \frac{b_5c_3}{b_3 n},
\end{eqnarray*}

\noindent the last inequality following from Lemma \ref{l:G_nk123}~(ii) (recalling that $Z_0 \le Y_0$). Consequently, for $1\le n\le M$,
$$
\E( m^{Z_n} \, {\bf 1}_{\{ Z_n\in [b_3 n, \, b_4 n]\}} \, {\bf 1}_{\{ \P( L_n^{(0)} \ge b_1 n^2 \, | \, Z_n) \ge b_5\}})
\ge
\frac{c_{37}}{n}
-
\frac{b_5c_3}{b_3 n}
=
\frac{c_{37}}{2n} \, ,
$$

\noindent if we choose $b_5 := \frac{c_{37}b_3}{2c_3}$. Let $B_n := \{y\in [b_3 n, \, b_4n] \cap \z: \, \P(Z_n=y)>0, \, \P(L_n^{(0)} \ge b_1 n^2 \, | \, Z_n=y)\ge b_5\}$. We have just proved that $\E( m^{Z_n} \, {\bf 1}_{\{Z_n \in B_n\}}) \ge \frac{c_{37}}{2n}$ for $1\le n\le M$. This yields Lemma \ref{l:B_nN_n>M^2} by taking $n=M$.\qed

\bigskip

\noindent {\it Proof of Lemma \ref{l:gen_fct_W}.} (i) We first prove that there exist constants $c_{33} \in (0, 1)$ and $c_{38}\in (0, \, 1)$, such that for any $s\in [m, \, m+\frac{c_{33}}{M}]$,
\begin{equation}
    \E(s^{W_0})-(s-1)\E(W_0s^{W_0})
    \ge 
    c_{38}.
    \label{pf_l:key}
\end{equation}

Let $s\ge m$. Since $W_0=Z_M \, {\bf 1}_{\{ Z_M \notin B_M\} }$ and $y\ge c_{29}\, M$ for any $y\in B_M$,
\begin{eqnarray*}
 && [(s-1)\E(Z_Ms^{Z_M})-\E(s^{Z_M})]-[(s-1)\E(W_0s^{W_0})-\E(s^{W_0})]
    \\
 &=& \E\{ [((s-1)Z_M-1)s^{Z_M}+1] \, {\bf 1}_{\{Z_M\in B_M\}} \}
    \\
 &\ge&  ((m-1)c_{29}M-1)\E(m^{Z_M} \, {\bf 1}_{\{ Z_M\in B_M\} }).
\end{eqnarray*}

\noindent We choose $M$ to be larger than $\frac{2}{(m-1)c_{29}}$, which ensures that $(m-1)c_{29}M-1\ge \frac12 (m-1)c_{29}M$. Furthermore, we choose $M$ to be sufficiently large so that we are entitled to apply Lemma~\ref{l:B_nN_n>M^2} to have $\E(m^{Z_M} \, {\bf 1}_{\{ Z_M\in B_M\} }) \ge \frac{c_{36}}{M}$. Combining these inequalities, we get, for $s\ge m$,
\begin{equation}
    [(s-1)\E(Z_Ms^{Z_M})-\E(s^{Z_M})]
    -
    [(s-1)\E(W_0s^{W_0})-\E(s^{W_0})]
    \ge
    \frac12(m-1)c_{29}c_{36}.
    \label{claim0}
\end{equation}

\noindent On the other hand, we claim the existence of constants $c_{39}>0$ and $c_{40}\ge 1$ such that for all $s\in [m, \, m+\frac{c_{39}}{M}]$,
\begin{equation}
    (s-1)\E(Z_M s^{Z_M})-\E(s^{Z_M})
    \le 
    c_{40}(s-m)M.
    \label{claim}
\end{equation}

To prove our claim, let $H_n(s) := \E(s^{Z_n})$ be the moment generating function of $Z_n$. We have, for $s\ge m$,
\begin{eqnarray*}
    (s-1) \E(Z_ns^{Z_n})
    -
    \E(s^{Z_n})
 &=& s(s-1)H_n'(s) - H_n(s)
    \\
 &\le& s(s-1) [H_n'(m)+(s-m)H_n''(s)] - H_n(m)
    \\
 &\le& (s-m)(s+m-1)H_n'(m) + s(s-1)(s-m)H_n''(s),
\end{eqnarray*}

\noindent where, in the last inequality, we have used the fact that $(m-1)\E(Z_nm^{Z_n}) - \E(m^{Z_n}) \le 0$. On the right-hand side, we have $H_n'(m)\le \frac{c_3}{m}$ (by Lemma \ref{l:G_nk123}~(ii) and the fact  that $Z_0 \le Y_0$), and use the existence of constants $c_{39}>0$ and $c_{41}\ge 1$ (see Theorem~4 and Corollary~1 of \cite{xz_stable}) such that $\E(Z_n^2 s^{Z_n}) \le c_{41}\, M$ for $1\le n\le M$ and $s\in [m, \, m + \frac{c_{39}}{M}]$; this implies $H_n''(s)\le c_{41} \, M$. Consequently, there exists a constant $c_{40}\ge 1$ 
such that for all $s\in [m, \, m+ \frac{c_{39}}{M}]$ and $1\le n\le M$,
\begin{equation}
    (s-1)\E(Z_ns^{Z_n})-\E(s^{Z_n})
    \le 
    c_{40}(s-m)M\, ,
    \label{pf_lemma5.5}
\end{equation}

\noindent which obviously implies the claim in \eqref{claim}.

We choose $c_{33} := \min\{\frac{1}{4c_{40}}(m-1)c_{29}c_{36}, \, c_{39}, \, \frac12\}$. If $s\in [m, \, m+\frac{c_{33}}{M}]$, then by \eqref{claim} and \eqref{claim0},
$$
(s-1)\E(W_0s^{W_0})-\E(s^{W_0})
\le
c_{40}(s-m)M -\frac{1}{2}(m-1)c_{29}c_{36} 
\le 
-\frac{1}{4}(m-1)c_{29}c_{36}.
$$

\noindent Setting $c_{38} := \min \{ \frac{1}{4}(m-1)c_{29}c_{36}, \, \frac12\}$ yields \eqref{pf_l:key}. 

(ii) We complete the proof of the lemma. Write $F_n(t)=\E(t^{W_n})$ for $t> 0$. Let $s:= m+\frac{c_{33}}{M}$. By \eqref{pf_l:key},
$$
F_0(s) - (s-1)sF_0'(s)
\ge 
c_{38}.
$$

\noindent By \eqref{p11}, $F_{n+1}(s) - (s-1)sF_{n+1}'(s) = [F_n(s) - m(s-1)F_n'(s)] F_n(s)^{m-1} \ge F_n(s) - s(s-1)F_n'(s)$, which implies that
$$
F_n(s)-(s-1)sF_n'(s)
\ge
F_0(s) - (s-1)sF_0'(s)
\ge 
c_{38},
\qquad
n\ge 0.
$$

\noindent Using $s(s-1) \ge m(m-1)$, this yields:
$$
F_n(s) 
\ge 
c_{38}
+
m(m-1) F_n'(s) ,
\qquad
n\ge 0.
$$

\noindent By convexity, $F_n'(s)\ge \frac{F_n(s)-F_n(m)}{s-m} = \frac{M}{c_{33}} (F_n(s)-F_n(m))$. Thus
$$
F_n(s) 
\ge 
c_{38}
+
c_{42} \, M (F_n(s)-F_n(m)) ,
\qquad
n\ge 0,
$$

\noindent with $c_{42} := \frac{m(m-1)}{c_{33}}$. Since $W_0\le Z_M\le Y_M$, we have $F_n(m)\le G_{n+M}(m)$, and the latter is bounded, according to Lemma \ref{l:G_nk123}~(iv), by $1+\frac{c_6}{n+M}$, a fortiori by $1+ \frac{c_6}{n}$ (for $n\ge 1$). Consequently,
$$
F_n(s)
\ge 
c_{38} 
+
c_{42} \, M \Big( F_n(s)-1-\frac{c_6}{n} \Big) ,
\qquad
n\ge 1\, ;
$$

\noindent in other words, we have
$$
F_n(s)
\le
\frac{c_{42} M(1+\frac{c_6}{n}) - c_{38}}{c_{42} M-1} 
=
1+ \frac{(1- c_{38}) + c_{42} M \frac{c_6}{n}}{c_{42} M-1} \, .
$$

\noindent We now take $n= n_0 := \lfloor c_{34} M \rfloor$, with $c_{34}:= \max\{ \frac{4c_6c_{42}}{c_{38}}, \, 2\}$. This choice ensures that $c_{42} M \frac{c_6}{n} \le c_{42} M \frac{c_6}{\frac12 \, c_{34} M} \le \frac12 c_{38}$, so $F_{n_0}(s) \le 1+ \frac{1- \frac12 c_{38}}{c_{42} M-1}$. This implies that
$$
\frac{m}{s}\, F_{n_0}(s)^{m-1}
=
\frac{F_{n_0}(s)^{m-1}}{1+ \frac{c_{33}}{mM}}
\le
\frac{1}{1+ \frac{c_{33}}{mM}} \Big( 1+\frac{1- \frac12 c_{38}}{c_{42} M -1} \Big)^{\! m-1} \, .
$$

\noindent Recall that $c_{42} = \frac{m(m-1)}{c_{33}}$. We choose $M_0$ so large that $\frac{1}{1+ \frac{c_{33}}{mM}} ( 1+\frac{1- \frac12 c_{38}}{c_{42} M -1})^{m-1} \le 1- \frac13 c_{38} \frac{c_{33}}{mM}$ for all $M\ge M_0$. Consequently, for $M\ge M_0$,
$$
\frac{m}{s}\, F_{n_0}(s)^{m-1}
\le 
1- \frac13 c_{38} \frac{c_{33}}{mM} \, ,
$$

\noindent proving Lemma \ref{l:gen_fct_W}.\qed

\subsection{Step 2. Study of $\P(Z_n \ge \ell+1)$}
\label{subs:P(Zn>k)}

The goal in this step is to prove Lemma \ref{l:P(Z>k)} below concerning $\P(Z_n \ge \ell+1)$. We start with an upper bound for the moment generating function of an arbitrary Derrida--Retaux system. The bound is simple, but not sharp in general. 

\medskip

\begin{lemma}
 \label{l:subcritcalspeed}
 
 Let $(X_n, \, n\ge 0)$ be a Derrida--Retaux system. Let $s>m$ and $\theta\in (0,1)$. If $\frac{m}{s} [\E(s^{X_0})]^{m-1}\le \theta$, then for $n\ge 0$,
 $$
 \E(s^{X_n}) \le 1 + (s-m)\theta^n.
 $$

\end{lemma}

\noindent {\it Proof.} Let $s>m$. Consider the sequence $(x_n, \, n\ge 0)$ defined by $x_0 := \E(s^{X_0})-1$ and $x_{i+1} := f(x_i)$ for $i\ge 0$, where $f(x) := \frac{(1+x)^m}{s}-\frac{1}{s}$ for $x\ge 0$. 

Since $f'(x)=\frac{m}{s}(1+x)^{m-1}$, we have, for all $x\in [0, \, x_0]$, $f'(x) \le f'(x_0) \le \theta$. By convexity, for $i\ge 0$, $x_{i+1} = f(x_i)-f(0) \le f'(x_i) \, x_i \le \theta x_i$. It follows that
\begin{equation}
    x_n \le \theta^n x_0,
    \qquad
    n\ge 0\, .
    \label{xn<}
\end{equation}

On the other hand, by definition of the Derrida--Retaux system, for $i\ge 0$, $\E(s^{X_{i+1}})=\frac{1}{s} [\E(s^{X_i})]^m + (1-\frac{1}{s}) [\P(X_i=0)]^m \le \frac{1}{s}[\E(s^{X_i})]^m + (1-\frac{1}{s})$, i.e., $\E(s^{X_{i+1}})-1\le f(\E(s^{X_i})-1)$. Iterating the inequality gives that
$$
\E(s^{X_n})-1\le x_n,
$$

\noindent which, according to \eqref{xn<}, is bounded by $\theta^n x_0$. By assumption, 
$$
\E(s^{X_0}) 
\le 
\Big(\frac{s\theta}{m}\Big)^{\! 1/(m-1)} 
\le 
\Big(\frac{s}{m}\Big)^{\! 1/(m-1)}
\le
\frac{s}{m} \, ,
$$

\noindent so $x_0 = \E(s^{X_0})-1 \le \frac{s}{m}-1 \le s-m$. Lemma \ref{l:subcritcalspeed} follows immediately.\qed

\bigskip

We have now all the ingredients to prove the main result for the Derrida--Retaux system $(Z_n, \, n\ge 0)$, stated as follows. 

\medskip

\begin{lemma}
 \label{l:P(Z>k)} 
 
 Assume $\E(Y_0^3m^{Y_0})<\infty$. There exist constants $c_{43}>0$, $c_{44}>0$ and sufficiently large $M_0$ 
 such that for all integers $n\ge M\ge M_0$, 
 \begin{equation}
     \P(Z_n \ge \ell+1)
     \le 
     c_{43} \Big( \frac{1}{M^2} + \frac{1}{M}\, \ee^{- c_{44}\frac{n}{M}} \Big) m^{-\ell},
     \qquad
     \ell\ge 0\, .
     \label{P(Z>k)}
 \end{equation}

\end{lemma}

\medskip

\noindent {\it Proof.} Let $M \ge M_0$, where $M_0$ is sufficiently large such that Lemma \ref{l:gen_fct_W} applies to $M\ge M_0$. Let $s = s(M) := m+\frac{c_{33}}{M}$ and $n_0 = n_0(M) := \lfloor c_{34}M \rfloor$. By Lemma \ref{l:gen_fct_W},
$$
\frac{m}{s} \, [\E(s^{W_{n_0}})]^{m-1}
\le 
1- \frac{c_{35}}{M}.
$$

\noindent This entitles us to apply Lemma \ref{l:subcritcalspeed} to $(W_{n+n_0}, \, n\ge 0)$, to see that
$$
\E(m^{W_{n+n_0}})-1
\le 
\E(s^{W_{n+n_0}})-1
\le
\frac{c_{33}}{M} (1- \frac{c_{35}}{M})^n,
\qquad
n\ge 0.
$$

Consider the system $(\widetilde{Z}_n := Z_{n+M}, \, n\ge 0)$. We define $\widetilde{L}_n^{(i)}$ the number of open paths at generation $n$ of this system with initial value $i$, exactly as $N_n^{(i)}$ for $(X_n)$ in \eqref{N}.\footnote{Notation: The letter $L$ is used for the number of open paths for $(Z_n)$, and $\widetilde{L}$ for $(\widetilde{Z}_n)$.} Let
$$
\widetilde{L}_n^{(B_M)}
:=
\sum_{i\in B_M} \widetilde{L}_n^{(i)} \, .
$$

\noindent Applying Lemma \ref{l:coupling} to $(\widetilde{Z}_k, \, k\ge 0)$ and $(W_k, \, k\ge 0)$, we get, for any $\ell\ge 1$,
$$
\P(Z_{n+M}=\ell, \; \widetilde{L}_n^{(B_M)} = 0)
\le 
\P(W_n=\ell)
\le 
\frac{\E(m^{W_n}-1)}{m^\ell-1}.
$$

\noindent Hence, for $n\ge 2n_0$, we have, for $\ell\ge 1$,
\begin{eqnarray*}
    \P(Z_{n+M}=\ell, \; \widetilde{L}_n^{(B_M)} = 0)
 &\le& \frac{1}{m^\ell-1}\, \frac{c_{33}}{M} \, (1- \frac{c_{35}}{M})^{n/2}
    \\
 &\le& \frac{1}{m^\ell-1}\, \frac{c_{33}}{M} \exp(-\frac{c_{35}n}{2M}).
\end{eqnarray*}

We now look at the situation $Z_{n+M}=\ell$ with $\widetilde{L}_n^{(B_M)} \ge 1$. On the event $\{ \widetilde{L}_n^{(B_M)} \ge 1\}$, we choose a vertex $v_M$ in generation $M$ of the system $(Z_n)$ such that $Z(v_m) \in B_M$ and that there is an open path from $v_M$ to generation $n+M$.\footnote{If there are several such vertices, we can take for example the leftmost one in the lexicographic order.} Call $L^{(0)} (v_M)$ the number of open paths in the system $(Z_n)$ from the initial generation to this particular vertex $v_M$ at generation $M$, with initial value $0$. Obviously,
\begin{eqnarray*}
    L_{n+M}^{(0)} \, {\bf 1}_{\{ Z_{n+M}=\ell\} }
 &\ge& L^{(0)} (v_M) \, {\bf 1}_{\{ \widetilde{L}_n^{(B_M)} \ge 1, \; Z_{n+M}=\ell\} } 
    \\
 &\ge& c_{31}M^2\, {\bf 1}_{\{ \widetilde{L}_n^{(B_M)} \ge 1, \; Z_{n+M}=\ell\} } \, {\bf 1}_{\{ L^{(0)} (v_M) \ge c_{31}M^2\} } \, .
\end{eqnarray*}

\noindent We take expectation on both sides. On the right-hand side, we first take conditional expectation with respect to the sigma-algebra $\mathscr{F}^Z_M$ generated by the values of the system $(Z_n)$ at generation $M$, and note that $\widetilde{L}_n^{(B_M)}$, $Z_{n+M}$, and $v_M$ (on the event $\{ \widetilde{L}_n^{(B_M)} \ge 1\}$) are all $\mathscr{F}^Z_M$-measurable. By definition of $B_M$, $\min_{y\in B_M} \P(L_M^{(0)} \ge c_{31}M^2 \, | \, Z_M=y) \ge c_{32}$. This leads to, for all $\ell\ge 1$,
\begin{eqnarray*}
    \E( L_{n+M}^{(0)} \, {\bf 1}_{\{ Z_{n+M}=\ell\} })
 &\ge& c_{31}M^2\, \E [ {\bf 1}_{\{ \widetilde{L}_n^{(B_M)} \ge 1, \; Z_{n+M}=\ell\} } \, \P (L^{(0)} (v_M) \ge c_{31}M^2 \, | \, \mathscr{F}^Z_M) ]
    \\
 &\ge& c_{31}c_{32}M^2 \, \P(\widetilde{L}_n^{(B_M)} \ge 1, \; Z_{n+M}=\ell) ,
\end{eqnarray*}

\noindent giving an upper bound for $\P( \widetilde{L}_n^{(B_M)} \ge 1, \; Z_{n+M}=\ell)$. Consequently, for $\ell\ge 1$,
\begin{eqnarray*}
    \P(Z_{n+M}=\ell)
 &=& \P( Z_{n+M}=\ell, \, \widetilde{L}_n^{(B_M)} \ge 1)
    +
    \P(Z_{n+M}=\ell, \, \widetilde{L}_n^{(B_M)} =0)
    \\
 &\le& \frac{\E(L_{n+M}^{(0)} \, {\bf 1}_{\{ Z_{n+M}=\ell\} })}{c_{31}c_{32}M^2}
    +
    \frac{1}{m^\ell-1}\, \frac{c_{33}}{M} \exp(-\frac{c_{35}n}{2M}).
\end{eqnarray*}

\noindent This yields that for $\ell\ge 1$,
$$
\P(Z_{n+M}\ge \ell)
\le
\frac{\E(L_{n+M}^{(0)} \, {\bf 1}_{\{ Z_{n+M}\ge \ell\} })}{c_{31}c_{32}M^2}
+
\frac{c_{33}}{M} \exp(-\frac{c_{35}n}{2M}) \sum_{j=\ell}^\infty \frac{1}{m^j-1} \, .
$$

\noindent Since $(Z_i, \, i\ge 0)$ is subcritical or critical, it follows from \eqref{E(N_Y>k)} that $\E(L_{n+M}^{(0)} \, {\bf 1}_{\{ Z_{n+M}\ge \ell\} }) \le \frac{c_1}{m^\ell}$; this implies the desired inequality \eqref{P(Z>k)} when $n\ge (2c_{34}+1)M$ (with $M\ge M_0$).

The case $M\le n< (2c_{34}+1)M$ (with $M\ge M_0$) is easy: we have $\E (m^{Z_n}\, {\bf 1}_{\{ Z_n \ge 1\} } ) \le \E (m^{Y_n}\, {\bf 1}_{\{ Y_n \ge 1\} } ) \le \frac{c_6}{n}$ (by Lemma \ref{l:G_nk123}~(iv)), thus for $\ell \ge 0$,
$$
\P(Z_n \ge \ell +1)
\le
\frac{\E (m^{Z_n}\, {\bf 1}_{\{ Z_n \ge 1\} })}{m^{\ell +1}} 
\le
\frac{c_6}{n}\, \frac{1}{m^{\ell +1}} 
\le
\frac{c_6}{M}\, \frac{1}{m^{\ell +1}} \, ,
$$

\noindent which implies again the desired inequality \eqref{P(Z>k)} because $n< (2c_{34}+1)M$. Lemma \ref{l:P(Z>k)} is proved.\qed

\subsection{Step 3. Proof of Proposition \ref{p:ub}}
\label{subs:pf_p_ub}

Assume $\E(Y_0^5m^{Y_0})<\infty$. Let $n$ and $M := \lfloor c_{44} \,\frac{n}{\log n}\rfloor$ (where $c_{44}>0$ is the constant in Lemma~\ref{l:P(Z>k)}) satisfy $n\ge M \ge M_0$, where $M_0$ is such that Lemma~\ref{l:P(Z>k)} applies to $n$ and $M$ and that $M\ge \frac{c_{44}}{2} \,\frac{n}{\log n}$. By Lemma \ref{l:P(Z>k)}, for $\ell\ge 0$,
\begin{eqnarray}
    \P(Z_n \ge \ell+1)
 &\le& c_{43} \Big( \frac{1}{M^2} + \frac{1}{M} \ee^{- c_{44}\frac{n}{M}} \Big) m^{-\ell}
    \nonumber
    \\
 &\le& c_{43} \Big( \frac{4(\log n)^2}{c_{44}^2 \, n^2} + \frac{2\log n}{c_{44}\, n} \, \ee^{-\log n} \Big) m^{-\ell}
    \nonumber
    \\
 &\le& \frac{c_{45}(\log n)^2}{n^2} \, \frac{1}{m^\ell} \, ,
    \label{pf_ub_eq1}
\end{eqnarray}

\noindent with $c_{45} := \frac{4c_{43}}{c_{44}^2} + \frac{2c_{43}}{c_{44}}$. 

Recall from \eqref{zeta} the definition: $\zeta(M) = -\log \E(Y_0^3m^{Y_0} \, {\bf 1}_{\{ Y_0>M \}})$. By the Markov inequality, $\E(Y_0^3m^{Y_0} \, {\bf 1}_{\{ Y_0>M \}}) \le \frac{c_{46}}{M^2}$, where $c_{46}:=\E(Y_0^5m^{Y_0})<\infty$. So $\zeta(M) \ge \log (\frac{M^2}{c_{46}})$ (and $\zeta(M)$ can be possibly infinite). Let $N_n^{(i)}$ be as before the number of open paths until generation $n$ in the system $(Y_n)$ with initial value $i$ (replacing $(X_n)$ by $(Y_n)$ in \eqref{N}) and let
$$
N_n^{(>M\zeta(M))} := \sum_{i> M \zeta(M)} N_n^{(i)} \, .
$$

\noindent By \eqref{recursion_E(YN)} and the Fubini--Tonelli theorem,
$$
\E((Y_n+1)m^{Y_n}N_n^{(>M\zeta(M))})
=
\E((Y_0+1)m^{Y_0}\, {\bf 1}_{\{Y_0>M \zeta(M)\}}) \prod_{i=0}^{n-1} G_i(m)^{m-1}.
$$

\noindent By Lemma \ref{l:G_nk123}~(v), $\prod_{i=0}^{n-1} G_i(m)^{m-1}\le c_8\, n^2$, whereas
$$
\E((Y_0+1)m^{Y_0}\, {\bf 1}_{\{Y_0>M \zeta(M)\}})
\le 
2 \, \E( Y_0 m^{Y_0}\, {\bf 1}_{\{Y_0>M \zeta(M)\}})
\le
\frac{2c_{46}}{(M \zeta(M))^4} \, ,
$$ 

\noindent with $c_{46}:=\E(Y_0^5m^{Y_0})$ as before. Hence
$$
\E((Y_n+1)m^{Y_n}N_n^{(>M\zeta(M))})
\le
\frac{2c_{46}}{(M \zeta(M))^4} \, c_8 \, n^2 \, .
$$ 

\noindent Therefore,
\begin{eqnarray}
    \P(Y_n \ge \ell+1, \, N_n^{(>M\zeta(M))}\ge 1) 
 &\le& m^{-(\ell+1)} \, \E((Y_n+1)m^{Y_n}N_n^{(>M\zeta(M))}) 
    \nonumber
    \\
 &\le& m^{-(\ell+1)} \, \frac{2c_{46}}{(M \zeta(M))^4} \, c_8 n^2 \, ,
    \label{pf_ub_eq2}
\end{eqnarray}

\noindent Applying Lemma \ref{l:coupling} to $(Y_k, \, k\ge 0)$ and $(Z_k, \, k\ge 0)$ yields that, for $\ell\ge 0$,
\begin{eqnarray*}
    \P(Y_n \ge \ell+1)
 &\le& \P(Y_n \ge \ell+1, \, N_n^{(>M\zeta(M))}\ge 1) 
    +
    \P(Z_n \ge \ell+1)
    \\
 &\le& m^{-(\ell+1)} \, \frac{2c_{46}}{(M \zeta(M))^4} \, c_8 n^2
    +
    \frac{c_{45}(\log n)^2}{n^2} \, \frac{1}{m^\ell} \, ,
\end{eqnarray*}

\noindent by means of \eqref{pf_ub_eq2} and \eqref{pf_ub_eq1}. Since $\zeta(M) \ge \log (\frac{M^2}{c_{46}})$ and $M := \lfloor \frac{c_{44}\, n}{\log n}\rfloor$, this implies Proposition \ref{p:ub} for sufficiently large $n$ (say  $n\ge n_0'$) and all $\ell \ge 0$. The case $2\le n<n_0'$ follows simply from the Markov inequality and Lemma \ref{l:G_nk123}~(i).\qed

\section{Proof of Theorem \ref{t:main}: lower bound}
\label{s:lb}

The lower bound in Theorem \ref{t:main} is a simple consequence of the following proposition.

\medskip

\begin{proposition}
 \label{p:lb}

 Assume $\E(t^{Y_0})<\infty$ for some $t>m$. Then for any $\lambda>0$, there exists a constant $\nu>0$ such that for all sufficiently large integer $n$,
 \begin{equation}
     \P(Y_n \ge \ell+1)
     \ge
     \frac{1}{n^2 (\log n)^\nu}\, \frac{1}{m^\ell} ,
     \qquad
     \forall \ell \in [0, \, \lambda n]\cap \z \, .
     \label{eq_p:lb}
 \end{equation}

\end{proposition}

\medskip

\noindent {\it Outline of the proof of Proposition \ref{p:lb}.} The technical ingredient in the proof is stated as Lemma \ref{l:E(S)>} below: there exist constants $\theta>0$ and $c_{48}>0$ such that for all sufficiently large integer $n$,
$$
 \sum_{k=0}^{k_n(\theta)-1} \mu_0(k) \E(S_n^{(k,\ell)})
 \ge
 \frac{c_{48}}{m^\ell} ,
 \qquad
 \forall \ell \in [0, \, \lambda n] \cap \z \, ,
 \leqno{\hbox{\it(\ref{eq_l:E(S)>})}}
$$
 
\noindent where
$$
 k_n(\theta) 
 :=
 \min \Big\{ i\ge 0: \, \E(Y_0^3 m^{Y_0} \, {\bf 1}_{\{ Y_0 \ge i+1\} }) \le \frac{\theta}{(\log n)^2} \Big\} ,
$$

\noindent and $S_n^{(k,\ell)}$ denotes as before (see \eqref{S}, with the system $(X_n)$ replaced by $(Y_n)$ there) the number of $(k, \, \ell)_n$--pivotal vertices in $\mathbb{T}_n$. 

Let us first describe how \eqref{eq_l:E(S)>} helps in the proof of Proposition \ref{p:lb}. Using a simple relation between pivotal vertices and open paths, we deduce from \eqref{eq_l:E(S)>} that for all sufficiently large $n$, there exist a constant $\varrho>0$ and a non-random integer $i^* = i^*(n) \in [0, \, c_{55} \log \log n]$ such that
$$
 \E(N_n^{(i^*)} \, {\bf 1}_{\{ Y_n \ge \ell+1\} })
 \ge
 \frac{1}{(\log n)^\varrho} \, \frac{1}{m^\ell} ,
 \qquad
 \forall \ell \in [0, \, \lambda n] \cap \z\, ,
 \leqno{\hbox{\it(\ref{eq_p:E(N)>})}}
$$

\noindent where, for any integer $i\ge 0$, $N_n^{(i)}$ denotes as before (see  \eqref{N}) the number of open paths in $\mathbb{T}_n$ to generation $n$ with initial value $i$, whereas $c_{55}>0$ is a suitable constant. [This is stated as Proposition \ref{p:E(N)>}, which is a weaker form of inequality \eqref{E(N_Y>k)} in the opposite direction.] 

Once \eqref{eq_p:E(N)>} is established, we proceed by coupling $(Y_i, \, i\ge 0)$ with some appropriate system $(X_i, \, i\ge 0)$ (depending on $n$) such that $Y_i \le X_i$ for all $i\ge 0$. A general inequality (Fact \ref{f:ub} below) borrowed from \cite{bmvxyz_conjecture_DR} for supercritical Derrida--Retaux systems says that
$$
\E(X_i) \le 3 ,
$$

\noindent for all $i$ smaller than a suitable constant multiple of $n$; moreover, a useful coupling inequality (called the ``bridge inequality"; see Fact \ref{f:bridge} below), also borrowed from \cite{bmvxyz_conjecture_DR}, says that 
$$
\E(N_n \, {\bf 1}_{\{ N_n \ge r\} } \, {\bf 1}_{\{ Y_n \ge \ell+1\} })
\le
\frac{2}{m^{k+\ell+1} \, \eta} \, \E(X_{n+k+\ell+1}) \, ,
$$

\noindent for all $n$, $k$, $\ell$ and $r$ lying in suitable domains, where $N_n$ denotes, as in \eqref{N}, the number of open paths in $\T_n$ to generation $n$ for $(Y_i)$. Adjusting the parameters appropriately, this leads to the following estimate: for any $a>0$ and $\lambda>0$, and all sufficiently large $n$,
$$
\E(N_n \, {\bf 1}_{\{ N_n \ge r\} } \, {\bf 1}_{\{ Y_n \ge \ell+1\} })
\le
\frac{1}{n^a} \, \frac{1}{m^\ell} ,
\qquad
\forall \ell \in [0, \, \lambda n] \cap \z\, ,
$$

\noindent where $r=r(n)$ is (approximately) a constant multiple of $n^2 \log n$. Together with \eqref{eq_p:E(N)>}, and using the trivial inequality $N_n \ge N_n^{(i^*)}$, this readily yields that
$$
\E(N_n^{(i^*)} \, {\bf 1}_{\{ N_n^{(i^*)} < r\} } \, {\bf 1}_{\{ Y_n \ge \ell+1\} })
\ge
\frac{1}{(\log n)^\varrho} \, \frac{1}{m^\ell}
-
\frac{1}{n^a} \, \frac{1}{m^\ell}
\ge
\frac12 \, \frac{1}{(\log n)^\varrho} \, \frac{1}{m^\ell} ,
$$

\noindent for all sufficiently large $n$, uniformly in $\ell \in [0, \, \lambda n] \cap \z$. Consequently, for all sufficiently large $n$ and uniformly in $\ell \in [0, \, \lambda n] \cap \z$,
$$
\P (Y_n \ge \ell+1)
\ge
\frac1r \, \E(N_n^{(i^*)} \, {\bf 1}_{\{ N_n^{(i^*)} < r\} } \, {\bf 1}_{\{ Y_n \ge \ell+1\} })
\ge
\frac{1}{2r} \, \frac{1}{(\log n)^\varrho} \, \frac{1}{m^\ell} ,
$$

\noindent from which Proposition \ref{p:lb} follows.

Let us now say some words about the proof of \eqref{eq_l:E(S)>}, which exploits the technique of pivotal vertices in Section \ref{s:pivot} and extends some technique to a more general context. We write
$$
\sum_{k=0}^{k_n(\theta)-1} \mu_0(k) \E(S_n^{(k,\ell)})
=
\sum_{k=0}^\infty \mu_0(k) \E(S_n^{(k,\ell)})
-
\sum_{k=k_n(\theta)}^\infty \mu_0(k) \E(S_n^{(k,\ell)}) \, .
$$

\noindent By Theorem \ref{t:new_iteration}, $\sum_{k=0}^\infty \mu_0(k) \E(S_n^{(k,\ell)}) = \mu_n(\ell)$, which is greater than or equal to $\frac{c_{47}}{m^\ell}$ for some constant $c_{47}>0$ (by Lemma \ref{l:truncated_moment_lb}; see \eqref{mu>} below). It remains to check that 
$$
\sum_{k=k_n(\theta)}^\infty \mu_0(k) \E(S_n^{(k,\ell)}) 
\le 
\frac{c_{47}}{2} \, \frac{1}{m^\ell} \, .
$$

\noindent We make a decomposition of the system at an intermediate generation $M = M(n) := \lfloor \frac{n}{1+c_{49}} \rfloor$, where $c_{49}>0$ is a (large) constant.\footnote{It is natural to introduce such a decomposition at an intermediate generation $M$, because the system at generation $M$ has ``nicer" properties than at the initial generation, which allows to use some renormalization idea starting from generation $M$. A similar approach was already used in the proof of the upper bound (Proposition \ref{p:ub}).} Contributions to $\E(S_n^{(k,\ell)})$ from the initial generation to generation $M$ are controlled (via the forthcoming Lemma \ref{l:E(N)>E(S)}) using the number of open paths to generation $M$, whereas contributions from generation $M$ to generation $n$ are investigated by means of Proposition \ref{p:ub}.

The rest of the section is devoted to the (rigorous) proof of Proposition \ref{p:lb}, presented in three steps. The first step is the technical part of the section; it gives precision on the number of pivotal vertices as in \eqref{eq_l:E(S)>} (Lemma \ref{l:E(S)>}). The second step gives a lower bound for $\max_{i\in [0, \, c_{55} \log \log n]} \E(N_n^{(i)} \, {\bf 1}_{\{ Y_n \ge \ell+1\} })$ as in \eqref{eq_p:E(N)>} (Proposition \ref{p:E(N)>}). The third step completes the proof of Proposition \ref{p:lb}.\qed 
 
\subsection{Step 1. Precision on the number of pivotal vertices}

Let $(Y_n, \, n\ge 0)$ be a critical Derrida--Retaux system: $\E(m^{Y_0}) = (m-1) \E(Y_0 m^{Y_0})<\infty$. Recall from \eqref{mu} that
\begin{eqnarray}
    \mu_n(k) 
 &:=& m^{-k}\, \E([1-(m-1)Y_n]m^{Y_n} \, {\bf 1}_{\{Y_n\le k\}})
    \nonumber
    \\
 &=& m^{-k}\, \E([(m-1)Y_n-1]m^{Y_n} \, {\bf 1}_{\{Y_n\ge k+1\}}) \, ,
    \label{mu_n=}
\end{eqnarray}

\noindent the equality \eqref{mu_n=} being a consequence of the fact that the system is critical, as we have already observed before. We have also noted (see \eqref{mu_bounds}) that 
$$
0\le \mu_n(k) \le m^{-k} ,
\qquad
k\ge 0, \; n\ge 0\, .
$$

\noindent On the other hand, by Lemma \ref{l:truncated_moment_lb}, if $\E(t^{Y_0}) <\infty$ for some $t>m$, then for any $\lambda>0$, there exists a constant $c_{47}>0$ such that for all sufficiently large integer $n$ and all $k \in [0, \, \lambda n]\cap \z$,
\begin{equation}
    \mu_n(k) \ge \frac{c_{47}}{m^k} \, .
    \label{mu>}
\end{equation}

\noindent Let $S_n^{(k,\ell)}$ denote, as in \eqref{S}, the number of $(k, \, \ell)_n$--pivotal vertices in $\mathbb{T}_n$. The aim of this step is to establish the following estimate.

\medskip

\begin{lemma}
\label{l:E(S)>}

 Assume $\E(t^{Y_0}) <\infty$ for some $t>m$. Let $\lambda>0$. There exist constants $\theta>0$ and $c_{48}>0$ such that for all sufficiently large integer $n$,
 \begin{equation}
     \sum_{k=0}^{k_n(\theta)-1} \mu_0(k) \E(S_n^{(k,\ell)})
     \ge
     \frac{c_{48}}{m^\ell} ,
     \qquad
     \forall \ell \in [0, \, \lambda n] \cap \z \, ,
     \label{eq_l:E(S)>}
 \end{equation}
 
 \noindent where
 $$
 k_n(\theta) 
 :=
 \min \Big\{ i\ge 0: \, \E(Y_0^3 m^{Y_0} \, {\bf 1}_{\{ Y_0 \ge i+1\} }) \le \frac{\theta}{(\log n)^2} \Big\} \, .
 $$

\end{lemma}

\medskip

The proof of the lemma relies on a simple preliminary result, stated as follows.

\medskip

\begin{lemma}
\label{l:E(N)>E(S)}

 For $n\ge 0$ and $i\ge k\ge 0$,
 $$
 \E [(Y_n+1)m^{Y_n} N_n^{(i)}]
 \ge
 m^{i-k} \, \P(Y_0=i) \sum_{\ell=0}^\infty (\ell+1)m^\ell \, \E(S_n^{(k,\ell)}) \, .
 $$

\end{lemma}

\medskip

\noindent {\it Proof of Lemma \ref{l:E(N)>E(S)}.} Let $n\ge 0$, $i\ge k\ge 0$ and $\ell \ge 0$ be integers. Consider $\sum_{v\in A_n^{(k,\ell)}} {\bf 1}_{\{ Y(v) = i\} }$. We have, by symmetry, $\E( \sum_{v\in A_n^{(k,\ell)}} {\bf 1}_{\{ Y(v) = i\} } ) = \P (Y_0 =i) \, \E(S_n^{(k,\ell)})$. On the other hand, $\sum_{v\in A_n^{(k,\ell)}} {\bf 1}_{\{ Y(v) = i\} } \le N_n^{(i)} \, {\bf 1}_{\{ Y_n = i-k+\ell\} }$. Hence
\begin{equation}
    \E(N_n^{(i)} \, {\bf 1}_{\{ Y_n = i-k+\ell\} })
    \ge
    \P (Y_0 =i) \, \E(S_n^{(k,\ell)}) ,
    \qquad
    i\ge k \ge 0\, .
    \label{E(N)>E(S)}
\end{equation}

\noindent Writing 
\begin{eqnarray*}
    \E((Y_n+1)m^{Y_n} N_n^{(i)}) 
 &=& \sum_{j=0}^\infty (j+1)m^j \, \E(N_n^{(i)} \, {\bf 1}_{\{ Y_n = j\} }) 
    \\
 &=& \sum_{\ell=-(i-k)}^\infty (i-k+\ell+1)m^{i-k+\ell}\, \E(N_n^{(i)} \, {\bf 1}_{\{ Y_n = i-k+\ell\} }) \, ,
\end{eqnarray*}

\noindent this yields that
\begin{eqnarray*}
    \E [ (Y_n+1)m^{Y_n} N_n^{(i)} ]
 &\ge& \sum_{\ell=-(i-k)}^\infty (i-k+\ell+1)m^{i-k+\ell}\, \P (Y_0 =i) \, \E(S_n^{(k,\ell)})
    \\
 &=& m^{i-k}\, \P (Y_0 =i) \sum_{\ell=-(i-k)}^\infty (i-k+\ell+1)m^\ell \, \E(S_n^{(k,\ell)}) \, .
\end{eqnarray*}

\noindent Since $i-k \ge 0$, this yields the desired inequality.\qed

\bigskip

We now have all the ingredients for the proof of Lemma \ref{l:E(S)>}.

\bigskip

\noindent {\it Proof of Lemma \ref{l:E(S)>}.} Fix $\lambda>0$. We write
$$
\sum_{k=0}^{k_n(\theta)-1} \mu_0(k) \E(S_n^{(k,\ell)})
=
\sum_{k=0}^\infty \mu_0(k) \E(S_n^{(k,\ell)})
-
\sum_{k=k_n(\theta)}^\infty \mu_0(k) \E(S_n^{(k,\ell)}) \, .
$$

\noindent On the right-hand side, $\sum_{k=0}^\infty \mu_0(k) \E(S_n^{(k,\ell)})$ equals $\mu_n(\ell)$ (by Theorem \ref{t:new_iteration}), which is greater than or equal to $\frac{c_{47}}{m^\ell}$ (by \eqref{mu>}) for all sufficiently large $n$ and uniformly in $\ell\in [0, \, \lambda n] \cap \z$. The proof of the lemma boils down to verifying that for some $\theta>0$ and all large $n$,
\begin{equation}
    \sum_{k=k_n(\theta)}^\infty \mu_0(k) \E(S_n^{(k,\ell)}) 
    \le 
    \frac{c_{47}}{2} \, \frac{1}{m^\ell},
    \qquad
    \forall \ell\in [0, \, \lambda n] \cap \z\, .
    \label{pf_l:E(S)>}
\end{equation}

Let $n>M\ge 0$ be integers. We extend the notion of pivotal vertices as follows: a vertex $v$ with $|v|=M$ is said to be $(k, \, \ell)_{M,n}$--pivotal if it is $(k, \, \ell)_{n-M}$--pivotal (in the sense of Definition \ref{def:pivot}) for the system $(Y_{M+i}, \, i\ge 0)$; in other words, $k+\xi_M(v) + \xi_{M+1}(v) +\cdots+\xi_{n-1}(v)=n-M+\ell$ and
$$
k+\xi_M(v) + \xi_{M+1}(v) + \cdots + \xi_i(v) \ge i-M+1, 
\qquad
\forall M\le i\le n-1,
$$

\noindent where, for $i\ge M$, $\xi_i(v) := \sum_{y\in \mathtt{bro}(v_i)} Y(y)$, with $\mathtt{bro}(v_i)$ denoting the set of the ``brothers" of $v_i$ (notation in agreement with \eqref{xi} if $M=0$); here, $v_i$ is the unique descendant of $v$ at generation $i$ (with $v_i :=v$ if $i=M$). We write $A_{M,n}^{(k, \ell)}$ for the set of all $(k, \, \ell)_{M,n}$--pivotal vertices in $\mathbb{T}_n$ (so $A_n^{(k, \ell)} = A_{0,n}^{(k, \ell)}$), and $S_{M,n}^{(k, \ell)} := | A_{M,n}^{(k, \ell)}|$. 
By definition,
$$
{\bf 1}_{\{ \mathfrak{e}_0 \in A_n^{(k,\ell)}\}}
=
\sum_{j=0}^\infty {\bf 1}_{\{ \mathfrak{e}_0 \in A_M^{(k,j)}\}} \, {\bf 1}_{\{ \mathfrak{e}_M\in A_{M,n}^{(j, \ell)}\}} \, .
$$

\noindent A similar idea has already been used in \eqref{A(n+1)_and_An}. As in Section \ref{s:pivot}, this implies that
$$
\E(S_n^{(k, \ell)})
=
m^n \, \E({\bf 1}_{\{ \mathfrak{e}_0 \in A_n^{(k,\ell)}\}})
=
m^n \sum_{j=0}^\infty \P ( \mathfrak{e}_0 \in A_M^{(k,j)}, \, \mathfrak{e}_M \in A_{M,n}^{(j, \ell)}) \, ;
$$

\noindent by independence of $\{ \mathfrak{e}_0 \in A_M^{(k,j)} \}$ and $\{ \mathfrak{e}_M \in A_{M,n}^{(j, \ell)} \}$, this leads to:
$$
\E(S_n^{(k, \ell)})
=
m^n \sum_{j=0}^\infty \P ( \mathfrak{e}_0 \in A_M^{(k,j)}) \, \P( \mathfrak{e}_M \in A_{M,n}^{(j, \ell)})
=
\sum_{j=0}^\infty \E(S_M^{(k,j)})\, \E(S_{M,n}^{(j, \ell)}) \, .
$$

\noindent Since $S_{M,n}^{(j, \ell)} =0$ for $j>n-M+\ell$, we obtain that 
\begin{equation}
    \E(S_n^{(k, \ell)})
    =
    \sum_{j=0}^{n-M+\ell} \E(S_M^{(k,j)})\, \E(S_{M,n}^{(j, \ell)}) \, .
    \label{pf:E(S)>_a}
\end{equation}

\noindent As such, the proof of \eqref{pf_l:E(S)>} will be complete if we are able to show that it is possible to choose an appropriate value of $M=M(n)$ such that for some $\theta>0$, all sufficiently large $n$, and all $\ell \in [0, \, \lambda n] \cap \z$,
\begin{eqnarray}
    \sum_{k=k_n(\theta)}^\infty \mu_0(k) \sum_{j=0}^{n-M+\ell-1} \E(S_M^{(k,j)})\, \E(S_{M,n}^{(j, \ell)})
    \le
    \frac{c_{47}}{4m^\ell} \, ,
    \label{pf:E(S)>_c}
    \\
    \sum_{k=k_n(\theta)}^\infty \mu_0(k) \E(S_M^{(k,j=n-M+\ell)})\, \E(S_{M,n}^{(j=n-M+\ell, \ell)})
    \le
    \frac{c_{47}}{4m^\ell} \, .
    \label{pf:E(S)>_b}
\end{eqnarray}

Let us first check \eqref{pf:E(S)>_b}. By the trivial inequalities $\sum_{k=k_n(\theta)}^\infty (\cdots) \le \sum_{k=0}^\infty (\cdots)$ and $S_{M,n}^{(n-M+\ell, \ell)} \le m^{n-M}$, we get
\begin{eqnarray*}
    \sum_{k=k_n(\theta)}^\infty \mu_0(k) \E(S_M^{(k,n-M+\ell)})\, \E(S_{M,n}^{(n-M+\ell, \ell)})
 &\le& m^{n-M} \sum_{k=0}^\infty \mu_0(k) \E(S_M^{(k,n-M+\ell)}) 
    \\
 &=& m^{n-M} \mu_M(n-M+\ell) \, ,
\end{eqnarray*}

\noindent the last equality being a consequence of Theorem \ref{t:new_iteration}. By Lemma \ref{l:G_nk123}~(iii), under the assumption $\E(Y_0^3m^{Y_0})<\infty$, we have $\E(Y_j^2 m^{Y_j})\le c_4\, j$ (for $j\ge 1$), so by the Markov inequality, there exists a constant $c_{49} >0$ such that for $j\ge 1$,
$$
\max_{i\ge c_{49}j} m^i \mu_j(i)
=
m^{\lceil c_{49}j\rceil} \mu_j(\lceil c_{49}j\rceil)
=
\E[((m-1)Y_j-1) m^{Y_j} \, {\bf 1}_{\{ Y_j\ge \lceil c_{49}j\rceil +1\} }]
\le
\frac{c_{47}}{4} \, .
$$

\noindent We choose
\begin{equation}
    M := \Big\lfloor \frac{n}{1+c_{49}} \Big\rfloor \, .
    \label{M}
\end{equation}

\noindent Then $m^{n-M+\ell} \mu_M(n-M+\ell) \le \frac{c_{47}}{4}$, proving \eqref{pf:E(S)>_b}. 

It remains to check \eqref{pf:E(S)>_c}. We bound $\E(S_{M,n}^{(j,\ell)})$. Let us write, for brevity, $\xi_i := \xi_i(\mathfrak{e}_0)$, $i\ge 0$. For $0\le j<n-M+\ell$,
\begin{eqnarray*}
    {\bf 1}_{\{ \mathfrak{e}_M\in A_{M,n}^{(j,\ell)} \} }
 &=& \sum_{i=0}^{j\wedge (n-M-1)} \; \sum_{r=1}^{n-M+\ell-j} {\bf 1}_{\{ \xi_M=0, \, \ldots, \xi_{M+i-1} =0, \, \xi_{M+i} = r\} } \, {\bf 1}_{\{ \mathfrak{e}_{M+i+1} \in A_{M+i+1,n}^{(j-i+r-1,\ell)}\} }
    \\
 &\le& \sum_{i=0}^{j\wedge (n-M-1)} \; \sum_{r=1}^{n-M+\ell-j} {\bf 1}_{\{ \xi_{M+i} = r\} } \, {\bf 1}_{\{ \mathfrak{e}_{M+i+1} \in A_{M+i+1,n}^{(j-i+r-1,\ell)}\} } \, .
\end{eqnarray*}

\noindent Hence
\begin{eqnarray*}
    \E(S_{M,n}^{(j,\ell)})
 &=& m^{n-M} \, \P(\mathfrak{e}_M\in A_{M,n}^{(j,\ell)})
    \\
 &\le& m^{n-M} \sum_{i=0}^{j\wedge (n-M-1)} \; \sum_{r=1}^{n-M+\ell-j} \P(\xi_{M+i} =r) \, \P(\mathfrak{e}_{M+i+1} \in A_{M+i+1,n}^{(j-i+r-1,\ell)})
    \\
 &=& \sum_{i=0}^{j\wedge (n-M-1)} \; \sum_{r=1}^{n-M+\ell-j} \P(\xi_{M+i} =r) \, m^{i+1} \, \E(S_{M+i+1,n}^{(j-i+r-1,\ell)}) \, .
\end{eqnarray*}

Let us have a closer look at $\P(\xi_{M+i} =r)$, or rather, at $\P(\xi_{M+i} \ge r)$. By definition, $\xi_{M+i}$ is the sum of $(m-1)$ independent copies of $Y_{M+i}$, and $Y_{M+i+1}$ is distributed as $(Y_{M+i}+\xi_{M+i}-1)^+$, with $Y_{M+i}$ and $\xi_{M+i}$ being independent. By Proposition \ref{p:ub}, $\P(Y_{M+i} \ge 1) \le c_{28} \, \frac{(\log (M+i))^2}{(M+i)^2}$, so 
$$
\P(\xi_{M+i} \ge 1) 
\le
(m-1) \, \P(Y_{M+i} \ge 1) 
\le 
(m-1)c_{28} \, \frac{(\log (M+i))^2}{(M+i)^2} \, .
$$

\noindent For $r\ge 2$, we use $\P(\xi_{M+i} \ge r) \le \frac{\P(Y_{M+i+1} \ge r-1)}{\P(Y_{M+i}=0)}$; since $\P(Y_{M+i}=0) \ge \frac12$ (for $M\ge M_0$ with a sufficiently large $M_0$), whereas $\P(Y_{M+i+1} \ge r-1) \le \frac{c_{28}}{m^{r-2}} \, \frac{(\log (M+i+1))^2}{(M+i+1)^2}$ (by Proposition \ref{p:ub} again), we get $\P(\xi_{M+i} \ge r) \le \frac{2c_{28}}{m^{r-2}} \, \frac{(\log (M+i+1))^2}{(M+i+1)^2}$. Since $M = M(n) := \lfloor \frac{n}{1+c_{49}} \rfloor$, we obtain, for some constant $c_{50} >0$ and all integer $r \ge 1$, 
$$
\P(\xi_{M+i} \ge r) \le \frac{c_{50}}{m^r} \, \frac{(\log n)^2}{n^2} \, .
$$

\noindent Hence, for $0\le j\le n-M+\ell-1$,
\begin{eqnarray*}
    \E(S_{M,n}^{(j,\ell)})
 &\le& \sum_{i=0}^{j\wedge (n-M-1)} \; \sum_{r=1}^{n-M+\ell-j} \frac{c_{50}}{m^r} \, \frac{(\log n)^2}{n^2} \, m^{i+1} \, \E(S_{M+i+1,n}^{(j-i+r-1,\ell)}) 
    \\
 &=& c_{50} \, \frac{(\log n)^2}{n^2}\, m^j \sum_{i=0}^{j\wedge (n-M-1)} \; \sum_{q=j-i}^{n-M+\ell-i-1} m^{-q}\, \E(S_{M+i+1,n}^{(q,\ell)}) \, ,
\end{eqnarray*}

\noindent with $q:= j-i+r-1$. Using $\sum_{q=j-i}^{n-M+\ell-i-1} (\cdots) \le \sum_{q=0}^{n-M+\ell} (\cdots)$, we arrive at: for $0\le j\le n-M+\ell-1$,
$$
\E(S_{M,n}^{(j,\ell)})
\le
c_{50} \, \frac{(\log n)^2}{n^2}\, m^j \sum_{i=0}^{j\wedge (n-M-1)} \;  \sum_{q=0}^{n-M+\ell} m^{-q}\, \E(S_{M+i+1,n}^{(q,\ell)}) \, .
$$

Recall from \eqref{mu>} that with our choice of $M$ as in \eqref{M},
$$
\mu_{M+i+1}(q)
\ge
c_{51} \, m^{-q} \, ,
$$

\noindent for some constant $c_{51}>0$, all sufficiently large $n$, all integers $0\le i\le n-M-1$ and $0\le q\le (1+\lambda)n$; in other words, $m^{-q} \le \frac{1}{c_{51}} \mu_{M+i+1} (q)$. This implies that 
\begin{eqnarray*}
    \E(S_{M,n}^{(j,\ell)})
 &\le& \frac{c_{50}}{c_{51}} \, \frac{(\log n)^2}{n^2}\, m^j \sum_{i=0}^{j\wedge (n-M-1)} \; \sum_{q=0}^{n-M+\ell} \mu_{M+i+1} (q) \, \E(S_{M+i+1,n}^{(q,\ell)}) 
    \\
 &\le& \frac{c_{50}}{c_{51}} \, \frac{(\log n)^2}{n^2}\, m^j \sum_{i=0}^j \sum_{q=0}^\infty \mu_{M+i+1} (q) \, \E(S_{M+i+1,n}^{(q,\ell)}) \, .
\end{eqnarray*}

\noindent By Theorem \ref{t:new_iteration} (applied to the system $(Y_{k+M+i+1}, \, k\ge 0)$ which is critical), we have $\sum_{q=0}^\infty \mu_{M+i+1} (q) \, \E(S_{M+i+1,n}^{(q,\ell)}) = \mu_n(\ell) \le \frac{1}{m^\ell}$. Hence, for all sufficiently large $n$, all $0\le \ell \le \lambda n$ and all $0\le j\le n-M+\ell -1$, with $c_{52} := \frac{c_{50}}{c_{51}}$, 
$$
\E(S_{M,n}^{(j,\ell)})
\le
c_{52} \, \frac{(\log n)^2}{n^2}\, \frac{(j+1) m^j}{m^\ell} \, .
$$

\noindent Since $\mu_0(k) := m^{-k} \, \E([(m-1)Y_0-1]m^{Y_0} \, {\bf 1}_{\{Y_0\ge k+1\}}) \le \frac{m-1}{m^k} \sum_{i = k+1}^\infty i m^i \, \P(Y_0=i)$, we get, for all sufficiently large $n$, all $k\ge 0$ and all $0\le \ell \le \lambda n$ (with $c_{53} := (m-1) c_{52}$)
\begin{eqnarray*}
 &&\mu_0(k) \sum_{j=0}^{n-M+\ell-1} \E(S_M^{(k,j)})\, \E(S_{M,n}^{(j, \ell)})
    \\
 &\le& \frac{c_{53}}{m^{k+\ell}} \frac{(\log n)^2}{n^2} \sum_{i = k+1}^\infty i m^i \, \P(Y_0=i) \sum_{j=0}^{n-M+\ell-1} (j+1) m^j \, \E(S_M^{(k,j)})
    \\
 &\le& \frac{c_{53}}{m^{k+\ell}} \frac{(\log n)^2}{n^2} \sum_{i = k+1}^\infty i m^i \, \P(Y_0=i) \sum_{j=0}^\infty (j+1) m^j \, \E(S_M^{(k,j)})\, .
\end{eqnarray*}

\noindent By Lemma \ref{l:E(N)>E(S)}, $\P(Y_0=i) \sum_{j=0}^\infty (j+1) m^j \, \E(S_M^{(k,j)}) \le m^{k-i} \, \E [(Y_M+1)m^{Y_M} N_M^{(i)}]$. Thus
$$
\mu_0(k) \sum_{j=0}^{n-M+\ell-1} \E(S_M^{(k,j)})\, \E(S_{M,n}^{(j, \ell)})
\le
\frac{c_{53}}{m^\ell}\, \frac{(\log n)^2}{n^2} \sum_{i = k+1}^\infty i \, \E [(Y_M+1)m^{Y_M} N_M^{(i)}] \, .
$$

\noindent Recall from \eqref{recursion_E(YN)} that $\E [(Y_M+1)m^{Y_M} N_M^{(i)}] = (i+1)m^i \, \P(Y_0=i) \prod_{j=0}^{M-1} G_j(m)^{m-1}$, where $G_n(s) := \E(s^{Y_n})$ as before. We have $\prod_{j=0}^{M-1} G_j(m)^{m-1} \le c_8 \, M^2 \le c_8 \, n^2$ (see Lemma \ref{l:G_nk123}~(v)). Therefore, with $c_{54} := c_{53} c_8$,
$$
\mu_0(k) \sum_{j=0}^{n-M+\ell-1} \E(S_M^{(k,j)})\, \E(S_{M,n}^{(j, \ell)})
\le
\frac{c_{54}}{m^\ell}\, (\log n)^2 \sum_{i = k+1}^\infty i (i+1)m^i \, \P(Y_0=i) \, .
$$

\noindent Summing over all integers $k\ge k_n(\theta)$ leads to:
\begin{eqnarray*}
 &&\sum_{k=k_n(\theta)}^\infty \mu_0(k) \sum_{j=0}^{n-M+\ell-1} \E(S_M^{(k,j)})\, \E(S_{M,n}^{(j, \ell)})
    \\
 &\le& \frac{c_{54}}{m^\ell}\, (\log n)^2 \sum_{k=k_n(\theta)}^\infty \sum_{i = k+1}^\infty i (i+1) m^i \, \P(Y_0=i)
    \\
 &=& \frac{c_{54}}{m^\ell}\, (\log n)^2 \sum_{i=k_n(\theta)+1}^\infty i (i+1) (i-k_n(\theta)) m^i \, \P(Y_0=i) \, ,
\end{eqnarray*}

\noindent the last equality being a consequence of the Fubini--Tonelli theorem. Since $i (i+1) (i-k_n(\theta)) \le i^2 (i+1) \le 2i^3$, and $\sum_{i=k_n(\theta)+1}^\infty 2i^3 m^i \, \P(Y_0=i) = 2\E( Y_0^3 m^{Y_0} \, {\bf 1}_{\{ Y_0 \ge k_n(\theta) +1\} })$, which is bounded by $\frac{2\theta}{(\log n)^2}$ by definition of $k_n(\theta)$, we see that, for all sufficiently large $n$ and all $0\le \ell \le \lambda n$,
$$
\sum_{k=k_n(\theta)}^\infty \mu_0(k) \sum_{j=0}^{n-M+\ell-1} \E(S_M^{(k,j)})\, \E(S_{M,n}^{(j, \ell)})
\le
\frac{2c_{54}\, \theta}{m^\ell} \, ,
$$

\noindent which is bounded by $\frac{c_{47}}{4m^\ell}$ if $\theta$ is chosen to satisfy $\theta \le \frac{c_{47}}{8c_{54}}$. This yields \eqref{pf:E(S)>_c}, and completes the proof of Lemma \ref{l:E(S)>}.\qed

\subsection{Step 2. A lower bound for $\max_i \E(N_n^{(i)} \, {\bf 1}_{\{ Y_n \ge \ell+1\} })$}

Let $(Y_n, \, n\ge 0)$ be a Derrida--Retaux system satisfying $\E(m^{Y_0}) = (m-1) \E(Y_0 m^{Y_0})<\infty$. For any $i\ge 0$, let $N_n^{(i)}$ denote as in \eqref{N} (associated with $(Y_n)$) the number of open paths in $\mathbb{T}_n$ to generation $n$ with initial value $i$. The aim of this step is to show that Lemma \ref{l:E(N)>E(S)} will give an interesting lower bound for $\E(N_n^{(i)})$ for a certain $i=i(n)$. 

\medskip

\begin{proposition}
\label{p:E(N)>}
                               
 Assume $\E(t^{Y_0})<\infty$ for some $t>m$. Let $\lambda>0$. There exist constants $\varrho>0$ and $c_{55}>0$ such that for all sufficiently large integer $n$, there exists a non-random integer $i^* = i^*(n) \in [0, \, c_{55} \log \log n]$ satisfying
 \begin{equation}
     \E(N_n^{(i^*)} \, {\bf 1}_{\{ Y_n \ge \ell+1\} })
     \ge
     \frac{1}{(\log n)^\varrho} \, \frac{1}{m^\ell} ,
     \qquad
     \forall \ell \in [0, \, \lambda n] \cap \z\, .
     \label{eq_p:E(N)>}
 \end{equation}

\end{proposition}

\medskip

\noindent {\it Proof.} By assumption, there exists $s_0>m$ such that $c_{56} := \E(Y_0^3 s_0^{Y_0})<\infty$. 

Let $\theta>0$ be the constant in Lemma \ref{l:E(S)>}. For all integer $k\ge 0$,
$$
\E(Y_0^3 m^{Y_0} \, {\bf 1}_{\{ Y_0 \ge k+1 \} })
\le
c_{56} \, \Big( \frac{m}{s_0} \Big)^{\! k+1} \, .
$$

\noindent We take $k := \lfloor c_{55} \log \log n \rfloor =: a_n$ (with $c_{55} := \frac{3}{\log (s_0/m)}$), so that $(\frac{m}{s_0})^{k+1} \le (\frac{m}{s_0})^{c_{55} \log \log n} = \frac{1}{(\log n)^3} \le \frac{\theta}{2c_{56} \, (\log n)^2}$ (for sufficiently large $n$). Hence for all sufficiently large $n$,
$$
\E(Y_0^3 m^{Y_0} \, {\bf 1}_{\{ Y_0 \ge a_n+1 \} })
\le
\frac{\theta}{2(\log n)^2} \, .
$$ 

\noindent On the other hand, $\E(Y_0^3 m^{Y_0} \, {\bf 1}_{\{ Y_0 \ge k_n(\theta) \} }) > \frac{\theta}{(\log n)^2}$ by definition. Thus for all sufficiently large $n$, we have $k_n(\theta) \le a_n$, and
$$
\E(Y_0^3 m^{Y_0} \, {\bf 1}_{\{ k_n(\theta) \le Y_0 \le a_n\} })
\ge
\frac{\theta}{(\log n)^2}
-
\frac{\theta}{2(\log n)^2}
=
\frac{\theta}{2(\log n)^2} \, .
$$

\noindent Since 
$$
\E(Y_0^3 m^{Y_0} \, {\bf 1}_{\{ k_n(\theta) \le Y_0 \le a_n\} }) 
=
\sum_{i=k_n(\theta)}^{a_n} i^3 m^i\, \P(Y_0=i) 
\le
a_n^4 m^{a_n} \max_{i\in [k_n(\theta), \, a_n] \cap \z} \P(Y_0=i) ,
$$

\noindent this yields that for all sufficiently large $n$, there exists $i^* = i^*(n) \in [k_n(\theta), \, a_n] \cap \z$, such that
$$
\P(Y_0=i^*)
\ge
\frac{\theta}{2(\log n)^2} \, \frac{1}{a_n^4 m^{a_n}} 
\ge
\frac{1}{(\log n)^{c_{57}}} \, ,
$$

\noindent with $c_{57} := 3+ c_{55} \log m$. 

On the other hand, by Lemma \ref{l:E(S)>}, for all sufficiently large $n$ and all $0\le \ell \le \lambda n$, there exists $k^* = k^*(n, \, \ell) \in [0, \, k_n(\theta)-1] \cap \z$, such that $\mu_0(k^*) \E(S_n^{(k^*,\ell)}) \ge \frac{c_{48}}{m^\ell k_n(\theta)}$; a fortiori $\E(S_n^{(k^*,\ell)}) \ge \frac{c_{48}}{m^\ell k_n(\theta)}$. Moreover, since $i^* \ge k_n(\theta) > k^*$, we are entitled to apply \eqref{E(N)>E(S)} to $i=i^*$ and $k=k^*$, to see that $\E(N_n^{(i^*)} \, {\bf 1}_{\{ Y_n = i^*-k^*+\ell\} }) \ge \P (Y_0 =i^*) \, \E(S_n^{(k^*,\ell)})$. Consequently, for all sufficiently large $n$,
$$
\E(N_n^{(i^*)} \, {\bf 1}_{\{ Y_n = i^*-k^*+\ell\} })
\ge
\P (Y_0 =i^*)\, \frac{c_{48}}{m^\ell k_n(\theta)}
\ge
\frac{1}{(\log n)^{c_{57}}}\, \frac{c_{48}}{m^\ell k_n(\theta)}
\ge
\frac{1}{(\log n)^{c_{57}}}\, \frac{c_{48}}{m^\ell a_n} \, .
$$

\noindent Since $i^*-k^* \ge 1$, we have $\E(N_n^{(i^*)} \, {\bf 1}_{\{ Y_n = i^*-k^*+\ell\} }) \le \E(N_n^{(i^*)} \, {\bf 1}_{\{ Y_n \ge \ell+1\} })$; hence for all sufficiently large $n$,
$$
\E(N_n^{(i^*)} \, {\bf 1}_{\{ Y_n \ge \ell+1\} })
\ge
\frac{1}{(\log n)^{c_{57}}}\, \frac{c_{48}}{m^\ell a_n}
=
\frac{1}{(\log n)^{c_{57}}}\, \frac{c_{48}}{m^\ell \lfloor c_{55} \log \log n \rfloor}\, ,
$$

\noindent proving the proposition. 
\qed

\subsection{Step 3. Proof of Proposition \ref{p:lb}}

The proof of Proposition \ref{p:lb} relies on two results borrowed from \cite{bmvxyz_conjecture_DR}, stated below as Facts \ref{f:bridge} and \ref{f:ub}. Fact \ref{f:bridge} is a bridge inequality connecting the number of open vertices in $(Y_n, \, n\ge 0)$ and the expected value of an appropriate supercritical system $(X_n, \, n\ge 0)$. Fact \ref{f:ub} gives an upper bound for $\E(X_n)$ when the system is ``slightly supercritical". 

\medskip

\begin{fact}
\label{f:bridge}

 {\bf (\cite[Theorem 4.2]{bmvxyz_conjecture_DR}; The bridge inequality)\footnote{The bridge inequality was stated in \cite{bmvxyz_conjecture_DR} for $N_n^{(0)}$ instead of $N_n$, but the proof was valid for $N_n$. Furthermore, $X_0$ and $Y_0$ were assumed in \cite{bmvxyz_conjecture_DR} to satisfy $\p(X_0=k) \ge \p(Y_0=k)$ for all integer $k\ge 1$ in order to ensure the so-called ``$XY$ coupling"; the latter holds automatically in our setting, thanks to the additional assumption $X_0 \ge Y_0$.}}
 Let $(X_n, \, n\ge 0)$ and $(Y_n, \, n\ge 0)$ be Derrida--Retaux systems such that $X_0 \ge Y_0$. Let $\eta >0$. If $\E( X_0-Y_0 \, | \, Y_0) \ge \eta$, then for all $r\ge 0$ and all integers $n\ge 0$, $\ell\ge 0$ and $k \in [0, \, \frac{r \eta}{2}]$,
 $$
 \E(X_{n+k+\ell})
 \ge
 \frac{m^{k+\ell} \, \eta}{2} \, \E \Big[ N_n \, {\bf 1}_{\{ N_n \ge r\} } \, {\bf 1}_{\{ Y_n \ge \ell\} } \Big] \, . 
 $$
\end{fact}

\medskip

\begin{fact}
 \label{f:ub}
 
 {\bf (\cite[Proposition 3.1 and Lemma 2.4]{bmvxyz_conjecture_DR})} 
 Assume $\E(m^{X_0})< (m-1)\E(X_0 \, m^{X_0}) <\infty$. Let $\delta_0 := (m-1) \E(X_0 \, m^{X_0}) - \E(m^{X_0})>0$. There exist constants $c_{58} \in (0, \, 1]$ and $c_{59}>0$, depending only on $m$, such that if $\delta_0 \in  (0, \, c_{58}]$, then
 $$
 \max_{0\le i\le c_{59} \, (\frac{\P (X_0=0)}{\delta_0})^{1/2}} \E (X_i) \le 3 \, .
 $$

\end{fact}

\medskip

\noindent {\it Proof of Proposition \ref{p:lb}.} Let $(Y_n, \, n\ge 0)$ be a critical system such that $\E(t^{Y_0})<\infty$ for some $t>m$. Let $\lambda>0$.

Let $n_2$ be such that Proposition \ref{p:E(N)>} applies to all $n\ge n_2$ (which is now taken for granted); in particular, Proposition \ref{p:E(N)>} yields, for all $n\ge n_2$, the existence of $i^* = i^*(n) \in [0, \, c_{55} \log \log n] \cap \z$ such that 
\begin{equation}
    \E(N_n^{(i^*)} \, {\bf 1}_{\{ Y_n \ge \ell +1\} }) 
    \ge 
    \frac{1}{(\log n)^\varrho} \, \frac{1}{m^\ell} \, ,
    \qquad
    \forall \ell \in [0, \, \lambda n] \cap \z \, .
    \label{i*}
\end{equation}

\noindent Let $\eta = \eta(n) := \frac{c_{60}}{n^2}$, where $c_{60} := \frac{c_{59}^2 \, \P(Y_0=0)}{2(\lambda+2)^2 c_{61}}$, with $c_{61} := m(m-1) \E(m^{Y_0})$. We can enlarge the value of $n_2$, if necessary, such that $\eta < \min\{ \frac12, \, \frac{c_{58}}{c_{61}}\}$ for all $n\ge n_2$. Let $U_0$ be a random variable independent of $Y_0$ and such that $\P(U_0 =1) = \eta = 1- \P(U_0=0)$. Let $X_0 := Y_0 + U_0$. Then $X_0 \ge Y_0$ and $\E( X_0-Y_0 \, | \, Y_0) = \eta$. So we are entitled to apply Fact \ref{f:bridge}. 

We also intend to apply Fact \ref{f:ub}, so let us note that by definition,
\begin{eqnarray*}
    \delta_0
 &:=&(m-1) \E(X_0 \, m^{X_0}) - \E(m^{X_0})
    \\
 &=& (m-1) \E(Y_0m^{Y_0}) \, \E(m^{U_0}) + (m-1) \E(m^{Y_0}) \, \E(U_0m^{U_0}) - \E(m^{Y_0}) \, \E(m^{U_0}) 
    \\
 &=& (m-1) \E(m^{Y_0}) \, \E(U_0m^{U_0})
    \\
 &=& c_{61} \, \eta \, ,  
\end{eqnarray*}

\noindent where $c_{61} := m(m-1) \E(m^{Y_0})$ as before. Moreover, $\P(X_0=0) = (1-\eta) \P(Y_0=0) \ge \frac12 \, \P(Y_0=0)$. By Fact \ref{f:ub},
$$
\max_{0\le i\le (\lambda+2)n}\E(X_i)
\le 
3\, .
$$

Fix $a>0$. Let $r = r(n):= \frac{2n^2}{c_{60} \log m} \, \log ( \frac{6n^{2+a}}{c_{60}})$ and $k = k(n):= \lfloor\frac{r \eta}{2}\rfloor$. Let $0\le \ell \le \lambda n$. Then $n+k+\ell+1 \le (\lambda+2)n$ (for $n\ge n_2$, with the value of $n_2$ enlarged if need be), so $\E(X_{n+k+\ell+1})\le 3$. On the other hand, by Fact \ref{f:bridge},
$$
\E(X_{n+k+\ell+1})
\ge
\frac{m^{k+\ell+1} \, \eta}{2}\, \E\Big[ N_n \, {\bf 1}_{\{ N_n \ge r\} } \, {\bf 1}_{\{ Y_n \ge \ell+1\} } \Big] \, ,
$$

\noindent which implies that
\begin{equation}
    \E\Big[ N_n \, {\bf 1}_{\{ N_n \ge r\} } \, {\bf 1}_{\{ Y_n \ge \ell+1\} } \Big]
    \le
    \frac{2}{m^{k+\ell+1} \, \eta}\, \E(X_{n+k+\ell+1})
    \le
    \frac{6}{m^{k+\ell+1} \, \eta} 
    \le
    \frac{1}{n^a}\, \frac{1}{m^\ell} \, ,
    \label{sept}
\end{equation}

\noindent recalling the definition of $k$, together with that of $r$ and $\eta$. We write
\begin{eqnarray*}
    \E\Big[ N_n^{(i^*)}\, {\bf 1}_{\{ N_n^{(i^*)} < r\} } \, {\bf 1}_{\{ Y_n \ge \ell+1\} } \Big]
 &=& \E [ N_n^{(i^*)}\, {\bf 1}_{\{ Y_n \ge \ell+1\} } ]
    -
    \E\Big[ N_n^{(i^*)}\, {\bf 1}_{\{ N_n^{(i^*)} \ge r\} } \, {\bf 1}_{\{ Y_n \ge \ell+1\} } \Big]
    \\
 &\ge& \E [ N_n^{(i^*)}\, {\bf 1}_{\{ Y_n \ge \ell+1\} } ]
    -
    \E\Big[ N_n\, {\bf 1}_{\{ N_n\ge r\} } \, {\bf 1}_{\{ Y_n \ge \ell+1\} } \Big] \, .
\end{eqnarray*}

\noindent By \eqref{i*}, $\E(N_n^{(i^*)} \, {\bf 1}_{\{ Y_n \ge \ell+1\} }) \ge \frac{1}{(\log n)^\varrho} \, \frac{1}{m^\ell}$, whereas by \eqref{sept}, $\E [ N_n \, {\bf 1}_{\{ N_n \ge r\} } \, {\bf 1}_{\{ Y_n \ge \ell+1\} } ] \le \frac{1}{n^a}\, \frac{1}{m^\ell}$. Therefore, for $n\ge n_2$ and $0\le \ell \le \lambda n$,
$$
\E\Big[ N_n^{(i^*)}\, {\bf 1}_{\{ N_n^{(i^*)} < r\} } \, {\bf 1}_{\{ Y_n \ge \ell+1\} } \Big]
\ge
\frac{1}{(\log n)^\varrho}\, \frac{1}{m^\ell}
- 
\frac{1}{n^a}\, \frac{1}{m^\ell}
\ge
\frac{1}{2(\log n)^\varrho}\, \frac{1}{m^\ell} \, .
$$

\noindent [We need to enlarge the value of $n_2$, if necessary, to ensure that $\frac{1}{n^a} \le \frac12 \, \frac{1}{(\log n)^\varrho}$ for all $n\ge n_2$.] On the other hand, $\E [ N_n^{(i^*)}\, {\bf 1}_{\{ N_n^{(i^*)} < r\} } \, {\bf 1}_{\{ Y_n \ge \ell+1\} } ] \le r \, \P(Y_n \ge \ell+1)$. Consequently, for $n\ge n_2$,
$$
\P(Y_n \ge \ell+1)
\ge
\frac1r \, \frac{1}{2(\log n)^\varrho}\, \frac{1}{m^\ell}
=
\frac{c_{60} \log m}{2n^2 \log ( \frac{6n^{2+a}}{c_{60}})}  \, \frac{1}{2(\log n)^\varrho}\, \frac{1}{m^\ell} \, ,
$$

\noindent from which \eqref{eq_p:lb} follows. Proposition \ref{p:lb} is proved.\qed

\bigskip
\bigskip

\noindent {\bf\Large Acknowledgements}

\bigskip

\noindent We are indebted to Victor Dagard, Bernard Derrida and Mikhail Lifshits for regular and stimulating discussions on the Derrida--Retaux model throughout the last few years. We also wish to thank two anonymous referees for their careful reading of the original manuscript and for their constructive comments.

\end{document}